\documentclass{article}

\usepackage{amssymb}
\def\suffix{ps}
\newcount\system
\global\system=3   % for unix(dvips)
\def\ifundefined#1{\expandafter\ifx\csname#1\endcsname\relax}
\ifundefined{figdir}\def\figdir{}\fi
\newcount\firstline
\newdimen\pswidth  \newdimen\xleft
\newdimen\psheight \newdimen\ytop \newdimen\ybot
\newcount\justx \newcount\justy
\global\justx=0 \global\justy=0
\newdimen\vpos \newtoks\labeL 
\newread\labeLfile \newdimen\xcoord \newdimen\ycoord
\newif\ifdoit 
\newbox\labox
\newdimen\xdvikwid 
\newdimen\xdvikht
\newdimen\pspoints
\newdimen\rwi
\newcount\temp
\def\readdim#1{\global\read\labeLfile to \temp
\global #1=\temp pt}
\def\figcrop#1{\par%  #1=filename
\openin\labeLfile=\figdir#1.lbl                                              
\global\read\labeLfile to\firstline\message{#1}               
\global\read\labeLfile to\temp%read overall dimensions                                     
\readdim{\ybot}
\readdim{\xleft}%               read upper left point
\readdim{\ytop}
\global\read\labeLfile to\justx%ignore
\global\read\labeLfile to\justy%ignore
\global\read\labeLfile to\labeL%ignore
\readdim{\pswidth}%            read lower right point
\global\advance\pswidth by -\xleft
\readdim{\psheight}
\global\advance\ybot by -\psheight
\global\advance\psheight by -\ytop
\global\read\labeLfile to\justx%ignore
\global\read\labeLfile to\justy%ignore
\global\read\labeLfile to\labeL%ignore                                    
\vbox to\psheight{\vfill
\ifnum\system=1
\fi %textures
\ifnum\system=2
\fi %msdos
\ifnum\system=3
                                                 \fi         %%unix:dvips
\ifnum\system=4
\fi         %%unix:dvips,scaled
\ifnum\system=1
\hbox to \pswidth{\kern-\xleft\special{postscriptfile \figdir#1.\suffix }\hfil}\fi
\ifnum\system=2
\hbox to \pswidth{\kern-\xleft\special{ps: plotfile \figdir#1.\suffix }\hfil}\fi
\ifnum\system=3
\hbox to \pswidth{\kern-\xleft\includegraphics{\figdir#1.\suffix}\hfil}\fi
\ifnum\system=4
\hbox to \pswidth{\kern-\xleft\includegraphics{\figdir#1.\suffix}\hfil}\fi
\ifnum\system=5
\hbox to \pswidth{\kern-\xleft\includegraphics{\figdir#1.\suffix}\hfil}\fi %orphee
\ifnum\system=6
   \xdvikwid=\pswidth
   \xdvikht=\psheight
   {\global\divide\xdvikwid by \pspoints}
   {\global\divide\xdvikht by \pspoints}
   \rwi=\xdvikwid
    {\global\multiply\rwi by 10}
\hbox to \pswidth{\kern-\xleft\includegraphics{\figdir#1.\suffix\space}\hfil}\fi                   %xdvik
\vskip -\baselineskip
\vskip -\ybot 
\vskip-\psheight %                                     
\hbox to\pswidth  {\hss}%                                            
\parindent=0pt\offinterlineskip                                       
\vpos=0 pt%                                                              
\loop\readdim{\xcoord}                                 
\ifdim \xcoord < -999pt \doitfalse\else\doittrue\fi                        
\ifdoit \advance \xcoord by -\xleft
\readdim{\ycoord}
\advance \ycoord by -\ytop                              
\global\read\labeLfile to\justx                                       
\global\read\labeLfile to\justy                                       
\global\read\labeLfile to\labeL
\global\setbox\labox=\hbox{\labeL\hskip-0.3em}%    
\advance\vpos by-\ycoord                                              
\vskip-\vpos \vpos=\ycoord                                         
\hbox to\pswidth{\hskip\xcoord %                                 
\hbox to 0pt{\ifnum\justx>0\hss\fi%                                   
\vbox to0pt{%                                                         
\ifnum\justy<2\vss\fi%                                                
\copy\labox\kern0pt%  
\ifnum\justy>0\vss\fi}%                                               
\ifnum\justx<2\hss\fi}%                                               
\hss}%                                                                
\repeat%                                                              
\advance\vpos by-\psheight%                                           
\vskip-\vpos %                                                     
}\closein\labeLfile}
\def\figplace#1#2#3{
\openin\labeLfile=\figdir#1.lbl
\ifeof \labeLfile
       \immediate\write16{***Can't find \figdir#1.lbl; Skipping it.***}
\else  \closein\labeLfile
       \null\hskip#2\raise #3 \hbox{\figcrop{#1}}
\fi
}
\def\figput#1{
\openin\labeLfile=\figdir#1.lbl
\ifeof \labeLfile
       \immediate\write16{***Can't find \figdir#1.lbl; Skipping it.***}
\else  \closein\labeLfile
       \hbox{\figcrop{#1}}
\fi
}

\def\NN{{\mathchoice {\mathrm{I \hspace{-0.2em} N}}
    {\mathrm{I \hspace{-0.2em} N}} {\mathrm{I \hspace{-0.14em} N}}
    {\mathrm{I \hspace{-0.14em} N}}}}

\def\QQ{{\mathchoice {\mathsf{I \hspace{-0.41em} Q}} {\mathsf{I
      \hspace{-0.47em} Q}} {\mathsf{I \hspace{-0.25em} Q}}
      {\mathsf{I \hspace{-0.2em} Q}}}}
\def\RR{{\mathchoice {\mathrm{I \hspace{-0.2em} R}} {\mathrm{I
        \hspace{-0.2em} R}} {\mathrm{I \hspace{-0.14em} R}}
    {\mathrm{I \hspace{-0.14em} R}}}}
\def\CC{\mathbb{C}}

\newfam\gotfam
\font\twlgot=eufm10 at 12pt \font\tengot=eufm10
 \font\sevengot=eufm7
\textfont\gotfam=\twlgot \scriptfont\gotfam=\tengot
\scriptscriptfont\gotfam=\sevengot

\newtheorem{prop}  {Proposition}
\newtheorem{lemma}  {Lemma}
\newtheorem{theor}   {Theorem}
\newtheorem{truc} {Definition}

\newcommand{\be}  {\begin{equation}}
\newcommand{\ee}  {\end{equation}}
\newcommand{\bea} {\begin{eqnarray}}
\newcommand{\eea} {\end{eqnarray}}
\newcommand{\lp}  {\left(}
\newcommand{\rp}  {\right)}

\newcommand{\Br}  {\overline}

\newcommand{\cB}  {{\cal B}}
\newcommand{\cL}  {{\cal L}}

\newcommand{\cC}  {{\cal C}}

\newcommand{\cO}  {{\cal O}}
\newcommand{\cM}  {{\cal M}}

\newcommand{\cN}  {{\cal N}}
\newcommand{\cT}  {{\cal T}}

\newcommand{\cF}  {{\cal F}}

\newcommand{\cA}  {{\cal A}}
\newcommand{\cE}  {{\cal E}}

\newcommand{\Om}  {\Omega}

\newcommand{\om}  {\omega}
\newcommand{\ta}  {\tau}

\newcommand{\si}  {\sigma}

\newcommand{\Ga}  {\Gamma}

\newcommand{\al}  {\alpha}

\newcommand{\rh}  {\rho}

\newcommand{\de}  {\delta}

\newcommand{\ph}  {\phi}

\newcommand{\GS}  {{\mathfrak S}}

\newcommand{\til}  {\tilde}

\def\Br{\overline}

\newcommand{\eqdef} {\stackrel{\rm def}{=}}

\def\endproof{\hfill\vrule height .6em width .6em depth
  0pt\goodbreak\vskip.25in}

\begin{document}

\title{Feynman Diagrams in Algebraic Combinatorics}

\author{Abdelmalek Abdesselam \\
\\
{\small LAGA, Institut Galil\'ee, CNRS UMR 7539}\\
{\small  Universit{\'e} Paris XIII}\\
{\small Avenue J.B. Cl{\'e}ment, F93430 Villetaneuse, France}\\
{\small email: abdessel@math.univ-paris13.fr}
}

\maketitle

%\hfill\eject

{\abstract{
We show, in great detail, how the perturbative tools of
quantum field theory allow one to rigorously obtain:
a ``categorified'' Faa di Bruno type formula for multiple composition,
an explicit formula for reversion and a proof of
Lagrange-Good inversion, all in the setting of multivariable
power series.
We took great pains to offer a self-contained presentation that,
we hope, will provide any mathematician who wishes,
an easy access to the wonderland of quantum field theory.
} }


%\hfill\eject

\section{Introduction}

In our last articles~\cite{Abdesselam1,Abdesselam2}, we showed the connection
between quantum field theory (QFT) and two research fields in pure
mathematics. The first field is the formal inverse approach
to the Jacobian conjecture which is a problem of commutative algebra;
and the second is the Lagrange-Good multivariable inversion which belongs to
enumerative combinatorics.
Although the arguments we succintly presented
in these articles
are mathematically rigorous,
it would be difficult to appreciate this fact without a sufficient mastery
of perturbation theory in QFT. 
There already exist a few carefully written introductions to QFT
aimed at a mathematical audience (\cite{Salmhofer}
is especially recommended for a general
introduction, as well as~\cite{FeldmanKT} for the particulars
of Fermionic theories, and since we are at the age of multimedia
one can also
{\em watch}~\cite{Brydges1}).
Nevertheless, we feel that these references
would benefit most the mathematical analyst,
rather than the algebraist or combinatorialist.
It is partly to fill this need, and also to develop an adequate mathematical
theory encompassing
our above mentioned papers, that this article was conceived.

Before we proceed, let us first explain what we mean by QFT and why
we feel that it is important for algebraic combinatorics.
As any theory in physics, QFT has both ``grammar'' and ``meaning''.
By ``grammar'' we refer to the mathematical structure of the theory,
independently of any physical interpretation.
The latter corresponds
to what we called the ``meaning''.
To give an analogy, the Navier-Stokes equation as an example of nonlinear
PDE, susceptible of a purely mathematical investigation, belongs to
the ``grammar'' of hydrodynamical theory; while
its physical interpretation, as describing the evolution of a real fluid
like air in the atmosphere, belongs to its ``meaning''.
Of course this distinction is an idealization, but it will help to
avoid misunderstandings in what follows.
Before we concentrate exclusively on the ``grammar'' of QFT let us
quickly though imperfectly define the ``meaning'' of QFT as the
description of the interaction via radiation fields
of constituents of matter,
considered as point-like objects, 
in a way that is compatible with the principles of quantum mechanics
and special relativity.
As such, it belongs to the rather specialized field of high
energy physics and it would be hard to justify its interest for
mathematicians at large.
However, much more is at stake concerning the ``grammar'' of QFT, which
we expect in the future to pervade
most fields of mathematics as it has those of physics (see~\cite{Witten}
for some prospective).
Indeed the ``grammar'' of QFT is fundamentally a {\em generalization of
calculus} and certainly the most exciting one
since Newton and Leibniz.

This generalization proceeds in two main directions:
functional or infinite-dimensional integration, and symbolic integration.
Roughly speaking, in the first direction one is interested
in defining ``natural'' measures on spaces of functions $\ph$ from
a {\em base} manifold $\cB$ to a {\em target} manifold
$\cT$. The calculus one learns in the first years of university
corresponds to the situation where $\cB$ is finite; the
integrals involved are the familiar ones in finitely many dimensions.
A one dimensional manifold $\cB$ corresponds for instance
to the Wiener measure and to stochastic processes related to Brownian
motion that are extensively studied in probability theory.
One truly starts doing QFT when the dimension of $\cB$ is at least two.
In fact, the most interesting and challenging situations occur in dimensions
2, 3 and 4, a pattern which surprisingly is also familiar
in topology.
The great difficulty of {\em constructive field theory}, which
is the branch of mathematical physics that adresses the problem of giving
a rigorous construction of these measures, 
and which has been honored by the choice of its most outstanding
problem among the prize problems of the Clay Foundation~\cite{JaffeW},
comes
from the requirement of ``naturality''.
The latter grosso modo means that the density of such measures, with
respect to the (ill-defined) Lebesgue measure (in case $\cT=\RR$ for
instance)
has to be defined only in terms of the {\em local} geometry
of $\cB$, $\cT$ and the maps $\ph:\cB\rightarrow\cT$ that one is
summing over.
For more on this we refer the reader to~\cite{Gawedzki,Salmhofer,Brydges2}.
We indulged in this digression because we are writing for
combinatorialists who might perhaps be agreeably surprised to learn that
the most successful methods to tackle this problem of mathematical analysis
are combinatorial!
Were it not already taken, a suitable denomination for constructive
field theory would be ``combinatorial analysis''.

We now come to the second direction of generalization we mentioned, that of
{\em symbolic integration}.
It will be the focus of this work which we hope will
deserve a place under the banner of~\cite{Cartier}.
A first and very fecund example of symbolic integral calculus
stemming from QFT is that of integration with
respect to anticommuting variables (Fermions).
It was introduced by F. A. Berezin, the rightful heir to Grassmann
and the elder Cartan.
One can define this operation, without speaking of ``integration'' at all,
as that of taking the ``top form'' coefficient in an element of
the exterior algebra of a finite-dimensional vector space.
Howerver one would then lose the {\em suggestive power}
of the integral notation among whose benefits has been the discovery of
 the celebrated Berezin change of
variable formula which underlies supersymmetry.
A ``bijective proof'' of this identity seems to us an urgent matter,
and an interesting question for the combinatorial community.
We will comment on this in section IV.
Besides, what physicists have discovered over the last
half-century are substantial fragments of a {\em dictionary} between
``integrals'' susceptible of a formal calculus (Feynman path-integrals),
and generating series in terms of discrete combinatorial structures
(Feynman diagrams).
The bridge between the two is given by Wick's theorem which we now state
in its {\em complex Bosonic} version.

\begin{theor}

Let $A\in\cM_n(\CC)$ be a matrix such that
$Re\ A\eqdef \frac{1}{2}(A+A^\ast)$ is positive definite.

\noindent
1) For any $J$ and $K$, two vectors in $\CC^n$,
one has
\be
\int_{\CC^n}d{\Br \ph}d{\ph}\ e^{-\ph^\ast A\ph+J^\ast\ph+\ph^\ast K}
=
\frac{e^{J^\ast A^{-1} K}}{det(A)}
\ee
where the $\ast$ means Hermitian conjugation,
$\ph\in\CC^n$ with components $\ph_1,\ldots,\ph_n$
is integrated with respect to the measure
\be
d{\Br \ph}d{\ph}
\eqdef
\prod_{i=1}^n
\frac{d(Re\ \ph_i)d(Im\ \ph_i)}{\pi}
\ee

\noindent
2) Let $i_1,\ldots,i_p$ and $j_1,\ldots,j_q$
be two collections of indices in $\{1,\ldots,n\}$, then
\be
\frac{\int_{\CC^n}d{\Br \ph}d{\ph}\ e^{-\ph^\ast A\ph}
\ph_{i_1}\ldots\ph_{i_p}{\Br\ph}_{j_1}\ldots{\Br\ph}_{j_q}}
{\int_{\CC^n}d{\Br \ph}d{\ph}\ e^{-\ph^\ast A\ph}}=
\left\{
\begin{array}{l}
0\ {\rm if}\ p\neq q \\
per\lp
(A_{i_\al j_\beta}^{-1})_{1\le\al,\beta\le p}
\rp\ {\rm if}\ p=q
\end{array}
\right.
\ee
where $per(M)$ denotes the permanent of the matrix $M$, and ${\Br\ph}_j$
is simply the complex conjugate of the component $\ph_j$ of $\ph$.

\noindent
3) For any polynomial in the $\ph_i$'s and the ${\Br\ph}_j$'s
considered as $2n$ independent variables, the effect
of integrating with respect to the Gaussian probability
measure
\[
\frac{d{\Br \ph}d{\ph}\ e^{-\ph^\ast A\ph}}
{\int_{\CC^n}d{\Br \ph}d{\ph}\ e^{-\ph^\ast A\ph}}
\]
on $\CC^n$ is the same as that of applying the
``differential operator''
\[
\exp\lp
\sum_{i,j=1}^n
\lp
\frac{\partial}{\partial\ph_i}
(A^{-1})_{ij}
\frac{\partial}{\partial{\Br\ph}_j}
\rp
\rp
\]
followed by the augmentation homomorphism, i.e. evaluation at
$\ph={\Br\ph}=0$.
\end{theor}

The proofs of {\it 2)} and {\it 3)} are easy consequences
of {\it 1)} which is an
exercise in ordinary calculus.
This theorem translates the evaluation of Gaussian integrals
into a combinatorial game, which should be evident from the
expansion of the permanent of {\it 2)} in terms of permutations or
by keeping track of which derivative acts on which factor
in the formulation {\it 3)}.
In what follows, we reverse the thrust and use the theorem as a
{\em definition} of ``integrals'' that can then be rigorously
constructed
in some rings of formal power series.
We are certainly not the first, and hopefully not the last,
to follow this line of
thought to make mathematical use of the
``folklore'' of perturbative QFT.
For instance, knot theorists have cashed in on this
idea~\cite{BarNatan1,BarNatan2},
and as we put this article into type,
we learned of~\cite{Fiorenza}, where
the principal
motivation is the study of the cohomology of moduli
spaces of curves, and where an effort similar to ours is made.
As for the present work, the most convenient mathematical
framework we found is the theory of combinatorial
species initiated by Joyal in his seminal work~\cite{Joyal}
(see also~\cite{BergeronLL}).
It provides us with the most reasonable compromise between
categorical ``abstract nonsense'' and the almost childlike simplicity
of Feynman diagrammatic
notation that has to be preserved at all costs.

Now let us outline the plan of this paper.
Part II is concerned with the formal calculus aspect
of the dictionary and addresses successively:
multiple composition of multivariable power series in II.1,
reversion in II.2,
and Lagrange-Good inversion with a, perhaps new, generalization of it
in II.3. The presentation in this part is deliberately heuristic.
Part III introduces a rigorous mathematical framework for the three
above mentioned topics which are treated
in the same sequence in III.1, III.2 and III.3 respectively.
We will end this article in section IV with a few comments indicating
some directions for further work.

\medskip
\noindent{\bf Aknowledgements :}
We thank D. Barsky, for introducing us to the combinatorial community
in particular to the S\'eminaire Lotharingien de Combinatoire,
and D.~Foata for encouraging us to write this detailed account
of the material we presented at the 49th SLC.
We thank A. Grigis for pointing out reference~\cite{GrigisS}.
We also thank J. Feldman for sharing his wonderful software for the drawing
of Feynman diagrams.
Finally the support of the Centre National de la Recherche Scientifique
is most gratefully aknowledged.

\section{Symbolic Integration}

\noindent{\bf Warning :} This section is heuristic.

\medskip
Let $R$ be a commutative ring with unit containing the field of rationals
$\QQ$.
Our concern in this section is to introduce a notion of integral calculus
for formal power series over $R$.
In the sequel a ``function''  $F:R^n\rightarrow R$ means
a formal power series $F\in R[[X_1,\ldots,X_n]]$.
Likewise a ``function'' $F:R^n\rightarrow R^n$ means
a system $F=(F_i)_{1\le i\le n}$ of $n$ power series in $R[[X_1,\ldots,X_n]]$.
If $u=(u_1,\ldots,u_n)$ is a vector of $n$ indeterminates, we
introduce the corresponding differential symbols
$du_1,\ldots,du_n$ and write $du\eqdef du_1\ldots du_n$ for their product.
Given two such vectors $u=(u_1,\ldots,u_n)$
and $v=(v_1,\ldots,v_n)$ we let $uv\eqdef u_1 v_1+\cdots+u_n v_n$.
We now introduce an integration symbol $\int$
in such a way that given a function $F:R^n\rightarrow R$,
$\int du F(u)$ is an element of $R$ we would like
to think of as the integral of $F$. We also introduce an
$n$-dimensional
``delta function''
$\de(u)=\de(u_1)\ldots\de(u_n)$, and postulate the following rules of
computation.

\medskip
\noindent{\bf Rule 1 :}
For any $F:R^n\rightarrow R$,
\be
\int du F(u) \de(u)= F(0)
\ee

\medskip
\noindent{\bf Rule 2 :}

\be
\int du\ e^{-uv}=\de(v)
\ee

\medskip
\noindent{\bf Rule 3 :}
All the rules of ordinary calculus are allowed:
integration by parts, Fubini, change of variables etc.
However, there is no absolute value for the Jacobian factor
in a change of variables.

\medskip
We will now use this apparently nonsensical calculational scheme for three
different applications.

\subsection{Composition}

Let $F$, $G$ be functions from $R^n$ to $R^n$, with no constant term

\medskip
\noindent{\bf Claim 1 :}
\be
(F \circ G)_i(X)=\int d{\Br s}ds
d{\Br t}dt d{\Br u}du
\ s_i
e^{-{\Br s}s-{\Br t}t-{\Br u}u
+{\Br s}F(t)+{\Br t}G(u)+{\Br u}X}
\label{FrondG}
\ee

\medskip
Here, $F=(F_i)_{1\le i\le n}$, $F_i(X)\in R[[X_1,\ldots,X_n]]$,
and likewise for $G$.
$F\circ G$ is the result of substituting $G_i(X)$ for $X_i$ in $F(X)$,
and $(F \circ G)_i(X)$
is the $i$-th component of the composition of $F$ with $G$.
Also ${\Br s}=({\Br s}_1,\ldots,{\Br s}_n)$,
$s=(s_1,\ldots,s_n)$,
${\Br t}=({\Br t}_1,\ldots,{\Br t}_n)$,
$t=(t_1,\ldots,t_n)$,
${\Br u}=({\Br u}_1,\ldots,{\Br u}_n)$
and
$u=(u_1,\ldots,u_n)$ are six vectors of indeterminates.
An expression like ${\Br s}s$ means $\sum_{i=1}^n {\Br s}_i s_i$,
and ${\Br s}F(t)\eqdef
{\Br s}_i F_i(t_1,\ldots,t_n)$.

The previous claim is obtained from the following symbolic calculation.
We first integrate over ${\Br u}$ in (\ref{FrondG})
according to Rule 2 :
\be
\int d{\Br u} e^{-{\Br u}u+{\Br u}X}
=\de(u-X)
\ee
therefore
\bea
\lefteqn{
\int d{\Br s}ds
d{\Br t}dt d{\Br u}du
\ s_i
e^{-{\Br s}s-{\Br t}t-{\Br u}u
+{\Br s}F(t)+{\Br t}G(u)+{\Br u}X}
} & & \nonumber \\
 & = & 
\int d{\Br s}ds
d{\Br t}dt du
\ s_i
e^{-{\Br s}s-{\Br t}t
+{\Br s}F(t)+{\Br t}G(u)}
\de(u-X)\\
 & = &
\int d{\Br s}ds
d{\Br t}dt dv
\ s_i
e^{-{\Br s}s-{\Br t}t
+{\Br s}F(t)+{\Br t}G(v+X)}
\de(v)
\eea
where we used the translation $v=u-X$, which gives by Rule 1
\[
\int d{\Br s}ds
d{\Br t}dt
\ s_i
e^{-{\Br s}s-{\Br t}t
+{\Br s}F(t)+{\Br t}G(X)}
\]
now this becomes after integration with respect to ${\Br t}$
\bea
\lefteqn{
\int d{\Br s}ds dt
\ s_i
e^{-{\Br s}s+{\Br s}F(t)}
\de(t-G(X))
} & & \nonumber \\
 & = & 
\int d{\Br s}ds dw
\ s_i
e^{-{\Br s}s+{\Br s}F(w+G(X))}
\de(w)
\eea
where $w=t-G(X)$, and finally we get
\bea
\int d{\Br s}ds
\ s_i
e^{-{\Br s}s+{\Br s}F(G(X))}
 & = & 
\int ds
\ s_i
\de(s-F(G(X))) \\
 & = & F_i(G(X))
\eea
Note that Rule 3 was used in the calculation by applying
``Fubini's theorem'' to perform the integrations in the chosen order,
and by using the change of variable formula in the simple case
of a translation.
It is easy to generalize this calculation to a multiple composition.
That is suppose $F^{(1)},\ldots,F^{(p)}$ are $p$
functions from $R^n$ to $R^n$ (with no constant term)
so that $F^{(j)}=(F_i^{(j)})_{1\le i\le n}$, with
each $F_i^{(j)}$ in $R[[X_1,\ldots,X_n]]$.
We can then write a similar formula for their $p$-fold
composition.

\medskip
\noindent{\bf Claim 2 :}
\bea
\lefteqn{
(F^{(1)}\circ\cdots\circ
F^{(p)})_i(X)=
\int \prod_{k=1}^{p+1}
(d{\Br u}^{(k)}du^{(k)})
} & & \nonumber\\
 & & u_i
\exp
\lp
-\sum_{k=1}^{p+1}
{{\Br u}^{(k)}u^{(k)}}
+\sum_{k=1}^p
{\Br u}^{(k)}F^{(k)}(u^{(k+1)})
+{\Br u}^{(p+1)}X
\rp
\eea
where, for each $k$, $1\le k\le p+1$,
${\Br u}^{(k)}=({\Br u}_1^{(k)},\ldots,{\Br u}_n^{(k)})$
and
$u^{(k)}=(u_1^{(k)},\ldots,u_n^{(k)})$ are vectors of 
indeterminates.
The above integral is formally over $2n(p+1)$ variables.

\subsection{Reversion}

Let $F$ be a function from $R^n$ to $R^n$
with no constant term and with invertible linear component.
We would like a formula for the compositional inverse
$F^{-1}:R^n\rightarrow R^n$.
For this purpose we let $\Om$ be any function $R^n\rightarrow R$
(not necessarily without constant term)
and we consider, with similar notation as in the previous section,
the integral
\[
\int d{\Br u}du\ \Om(u)
e^{-{\Br u}F(u)+{\Br u}Y}
\label{testfunc}
\]
We use Rule 3 to perform the change of variable
$v=F(u)-Y$
or $u=F^{-1}(v+Y)$.
If $F=(F_i(X))_{1\le i\le n}$ with $X=(X_1,\ldots,X_n)$
and $F_i(X)\in R[[X_1,\ldots,X_n]]$, for $1\le i\le n$,
we use the notation
\be
\partial F(Z)\eqdef
\lp
\frac{\partial F_i}{\partial X_j}(Z)
\rp_{1\le i,j\le n}
\ee
for the Jacobian matrix of $F$ ``at the point''
$Z=(Z_1,\ldots,Z_n)$.
It is an element of $\cM_n(R[[Z_1,\ldots,Z_n]])$.
The change of variable formula asserts that
in the above integral we can replace the dummy
integration variable $u$ by $F^{-1}(v+Y)$
and $du$ by
$det[\partial(F^{-1})(v+Y)] dv$.
Therefore
\bea
\int d{\Br u}du\ \Om(u)
e^{-{\Br u}F(u)+{\Br u}Y}
& = &
\int d{\Br u}dv\ det[\partial(F^{-1})(v+Y)] \nonumber\\
 & & \Om(F^{-1}(v+Y))
e^{-{\Br u}v} \\
 & = &
\int dv
\ det[\partial(F^{-1})(v+Y)]
\nonumber\\
 & & \Om(F^{-1}(v+Y))
\de(v)
\eea
by Rule 2.
Finally Rule 1 gives
\be
\int d{\Br u}du\ \Om(u)
e^{-{\Br u}F(u)+{\Br u}Y}
=
\Om(F^{-1}(Y))det[\partial(F^{-1})(Y)]
\ee
Applying this last formula successively to $\Om(u)\eqdef 1$
and $\Om(u)\eqdef u_i$, and noting the cancellation of the
determinantal factor, we obtain

\medskip
\noindent{\bf Claim 3 :}

\be
F^{-1}(Y)_i=
\frac{\int d{\Br u}du\ u_i
e^{-{\Br u}F(u)+{\Br u}Y}}
{\int d{\Br u}du\ e^{-{\Br u}F(u)+{\Br u}Y}}
\ee
which is a formula for the $i$-th component of the compositional inverse
$F^{-1}$ applied to $Y=(Y_1,\ldots,Y_n)$.
It is a generalization, due to V. Rivasseau
and the author~\cite{Abdesselam1},
of a formula that first appeared in the context of KAM
theory~\cite{Gallavotti}.
Note the compelling probabilistic interpretation as the average of $u_i$
with respect to the probability measure
proportional to
$d{\Br u}du\ e^{-{\Br u}F(u)+{\Br u}Y}$.

\subsection{Lagrange-Good inversion}

Let $G=(G_i)_{1\le i\le n}$ be a given system of formal power series
in $n$ indeterminates, and
let $F=(F_i)_{1\le i\le n}$ be the system of formal power series in
$R[[X_1,\ldots,X_n]]$, without constant term, implicitely defined
by the equations
\be
F_i=X_i G_i(F)\ \ {\rm for}\ \ 1\le i\le n
\ee
The implicit form of the multivariable Lagrange-Good
inversion formula says that for any $\Om:R^n\rightarrow R$
the coefficient
of $\frac{X^M}{M!}$
in
\[
\Om(F)\frac{1}{det(\de_{ij}-X_i\partial_j G_i(F))}
\]
is equal to the coefficient of $\frac{X^M}{M!}$
in $\Om(u)G(u)^M$.
Here we used the multiindex notation $M\eqdef(M_1,\ldots,M_n)\in\NN^n$,
$u^M\eqdef u_1^{M_1}\ldots u_n^{M_n}$,
$M!\eqdef M_1!\ldots M_n!$
and
$G(u)^M\eqdef G_1(u)^{M_1}\ldots G_n(u)^{M_n}$.
If $G_i=G_i(u)$,
is originally given in terms of the $u$ variables,
$\partial_j G_i(F)$ denotes the substitution of $u$ by $F=F(X)$
in $\frac{\partial G_i}{\partial u_j}(u)$.
Therefore
\be
(\de_{ij}-X_i\partial_j G_i(F))_{1\le i,j\le n}\in \cM_n(R[[X_1,\ldots,X_n]])
\ee
The Lagrange-Good formula can be compactly written as
\be
\Om(F)\frac{1}{det(\de_{ij}-X_i\partial_j G_i(F))}=
\sum_{M\in \NN^n}\frac{X^M}{M!}
\left.
\lp\frac{\partial}{\partial u}\rp^M
\right|_{u=0}
\Om(u)G(u)^M
\label{LG}
\ee

We will derive this identity using our symbolic calculus.
We consider as before the integral
\[
\int d{\Br u}du\ \Om(u)
e^{-{\Br u}u+{\Br u}XG(u)}
\]
Here
${\Br u}XG(u)\eqdef\sum_{i=1}^n
{\Br u}_i X_i G_i(u)$. In order to be able to apply Rule 2
to integrate over $\Br u$, we need first to perform the change of variables
$v=H(u)$ where $H_i(u)=u_i-X_i G_i(u)$. $H$ is considered as a function of
$u$ from $R^n$ to $R^n$.
$X$ plays the role of a parameter. Note that $F(X)=H^{-1}(0$).
Now $v=H(u)$ implies $u=H^{-1}(v)$ and
$du=det[\partial(H^{-1})(v)] dv$, therefore by Rule 3
\bea
\int d{\Br u}du\ \Om(u)
e^{-{\Br u}u+{\Br u}XG(u)}
 & = &\int d{\Br u}du\ \Om(u)
e^{-{\Br u}H(u)} \\
 & = & \int d{\Br u}dv\ \Om(H^{-1}(v))
det[\partial(H^{-1})(v)]
e^{-{\Br u}v}
\eea
Now by Rules 2 and 1 this becomes
\be
\int dv\ \Om(H^{-1}(v))
det[\partial(H^{-1})(v)]
\de(v)
= 
\Om(H^{-1}(0))
det[\partial(H^{-1})(0)]
\ee
\be
=  \Om(H^{-1}(0))
\frac{1}{det[\partial H(H^{-1}(0))]}
\ee
\be
=  \Om(F)\frac{1}{det(I-X\partial G(F))}
\ee
where $I$ is the $n\times n$ identity matrix and
$[X\partial G(F)]_{ij}\eqdef X_i \partial_j G_i(F)$.
Therefore
\be
\int d{\Br u}du\ \Om(u)
e^{-{\Br u}u+{\Br u}XG(u)}
=
\Om(F)\frac{1}{det(\de_{ij}-X_i\partial_j G_i(F))}
\label{leftLG}
\ee

There is a second way to do the computation, namely to expand the
exponential.
In multiindex notation it reads
\be
\int d{\Br u}du\ \Om(u)
e^{-{\Br u}u+{\Br u}XG(u)}
=
\sum_{M\in\NN^n}
\frac{X^M}{M!}
\int d{\Br u}du\ {\Br u}^M e^{-{\Br u}u}
\Om(u) G(u)^M
\ee
Now note that ${\Br u}^M e^{-{\Br u}u}={(-\frac{\partial}
{\partial u})}^M e^{-{\Br u}u}$
therefore one can integrate by parts
\bea
\int d{\Br u}du\ {\Br u}^M e^{-{\Br u}u}
\Om(u) G(u)^M
& = &
\int d{\Br u}du\ e^{-{\Br u}u}
\lp\frac{\partial}
{\partial u}\rp^M
\left[\Om(u) G(u)^M\right]\\
 & = & \int du\ \de(u)
\lp\frac{\partial}
{\partial u}\rp^M
\left[\Om(u) G(u)^M\right]\\
 & = &
\left.\lp\frac{\partial}
{\partial u}\rp^M\right|_{u=0}
\left[\Om(u) G(u)^M\right]
\eea
therefore
\be
\int d{\Br u}du\ \Om(u)
e^{-{\Br u}u+{\Br u}XG(u)}
=
\sum_{M\in\NN^n}
\frac{X^M}{M!}
\left.\lp\frac{\partial}
{\partial u}\rp^M\right|_{u=0}
\left[\Om(u) G(u)^M\right]
\label{rightLG}
\ee
Now (\ref{LG}) follows from (\ref{leftLG}) and (\ref{rightLG}).
Note also that by applying (\ref{leftLG}) successively to $\Om(u)\eqdef 1$
and $\Om(u)\eqdef u_i$, and computing the ratio we obtain a formula for
$F_i(X)$ which is

\medskip
\noindent{\bf Claim 4 :}
\be
F_i(X)=
\frac{\int d{\Br u}du\ u_i
e^{-{\Br u}u+{\Br u}XG(u)}}
{\int d{\Br u}du\ e^{-{\Br u}u+{\Br u}XG(u)}}
\ee

\medskip
Note the similarity with Claim 3.
In fact we derived Claim 4 from Claim 3 in~\cite{Abdesselam2}.
Note also that one of the advantages of our symbolic calculus is
that it can suggest generalizations and variations on the Lagrange-Good
formula. For instance, given $G=(G_i)_{1\le i \le n}$,
with $G_i=G_i(u)$, $u=(u_1,\ldots,u_n)$, and given
$n^2$ indeterminates
$(X_{ij})_{1\le i,j\le n}$, one can show that there is a unique solution
$F=(F_i)_{1\le i \le n}$, with $F_i\in R[[(X_{ij})_{1\le i,j\le n}]]$,
without constant term, to the equations
\be
F_i=\sum_{j=1}^n X_{ij} G_j(F)\ \ {\rm for}\ \ 1\le i\le n
\ee
The solution is given by the same expression as in claim 4, except that
\be
{\Br u}X G(u)\eqdef \sum_{i,j=1}^n {\Br u}_i
X_{ij} G_j(u_1,\ldots,u_n)
\ee
One also has a Lagrange-Good type identity which seems to be new, and
is our first encounter here with the inadequacy of the multiindex notation.
With $\Om:R^n\rightarrow R$, it reads
\[
\Om(F)\frac{1}{det(I-X\partial G(F))}=\sum_{k\ge 0}
\sum_{{i_1,\ldots,i_k,}\atop{j_1,\ldots,j_k=1}}^n
\frac{X_{i_1 j_1}\ldots X_{i_k j_k}}{k!}
\]
\be
\int d{\Br u}du\ \Om(u)
e^{-{\Br u}u}
{\Br u}_{i_1}\ldots{\Br u}_{i_k}
G_{j_1}(u)\ldots G_{j_k}(u)
\ee
\be
=  \sum_{k\ge 0}
\sum_{{i_1,\ldots,i_k,}\atop{j_1,\ldots,j_k=1}}^n
\frac{X_{i_1 j_1}\ldots X_{i_k j_k}}{k!}
\left.
\frac{\partial}{\partial u_{i_1}}
\ldots
\frac{\partial}{\partial u_{i_k}}
\right|_{u=0}
\left[
\Om(u)
G_{j_1}(u)\ldots G_{j_k}(u)
\right] 
\ee

\section{Feynman diagrams}

In contrast to the previous heuristic but conceptually important section,
we will now do some mathematics.
Throughout the remainder of this article $[n]$ will denote the 
set of the first $n$ nonnegative integers; and $\#(E)$ will denote
the cardinal of a finite set $E$.
As we will constantly use the notion of summable families
in power series rings, the reader who needs it should
consult~\cite{Bourbaki}
for a refresher.
Let $R$ be a commutative ring with unit containing  $\QQ$.
A single power series $F\in R[[X_1,\ldots,X_n]]$ is usually
specified using multiindex notation as
\be
F=\sum_{\al\in\NN^n}\frac{1}{\al!}u_\al X^\al
\label{mindex}
\ee
where $\al=(\al_1,\ldots,\al_n)$
is a multiindex in $\NN^n$, $\al!=\al_1!\ldots\al_n!$,
$X^\al=X_1^{\al_1}\ldots X_n^{\al_n}$, and $(u_\al)_{\al\in\NN^n}$
is the family of coefficients
defining $F$.
Note the normalization by $\frac{1}{\al!}$
which might seem like an insignificant matter of
convention but will be a leitmotiv in the following exposition,
namely to normalize by the cardinal of the group of {\em ambiguity}.
We will make this precise later.
Multiindex notation is more than highly impractical for our purposes,
it is in fact an example of a ``bad decategorification'' in the sense
of~\cite{BaezD} who were also influenced
by Joyal's theory of species.
It is rather more natural to use a {\em tensorial} notation
for our power series $F$ as
\be
F=\sum_{d\ge 0}\frac{1}{d!}
\sum_{i_1,\ldots,i_d}^n
F^{[d]}_{i_1\ldots i_d}
X_{i_1}\ldots X_{i_d}
\label{tensor}
\ee
Note again the $\frac{1}{d!}$
as our ambiguity group here is the symmetric group $\GS_d$.
If $(i_1,\ldots,i_d)\in [n]^d$
let $\mu(i_1,\ldots,i_d)\in\NN^n$
be the associated multiplicity multiindex that is
$\mu(i_1,\ldots,i_d)=(\mu_1,\ldots,\mu_n)$
where, for each $i$, $1\le i\le n$, $\mu_i$ is the number
of indices $r\in [d]$ such that $i_r=i$.
The translation between (\ref{mindex}) and (\ref{tensor}) is of course
\be
F^{[d]}_{i_1\ldots i_d}=
u_{\mu(i_1,\ldots,i_d)}
\ee
for any $d\ge 0$ and $(i_1,\ldots,i_d)\in [n]^d$.
$F^{[d]}_{i_1\ldots i_d}$ can be thought of as a tensor element
(the multidimensional analog of a matrix entry)
 of a symmetric $d$-covariant tensor (i.e. an element of the 
$R$-module $Sym^d((R^n)^\ast)$).
Feynman diagrams are first of all, the most efficient notation for
tensors and tensor contraction (see the Appendix of~\cite{PenroseR}
and also~\cite{Cvitanovic}) and the most
natural step after Einstein's convention of summing over repeated
indices.
We will for instance write
\be
\figplace{dessin1}{0 in}{-0.45 in}
\eqdef F^{[d]}_{i_1\ldots i_d}
\ee
\be
\figplace{dessin2}{0 in}{-0.47 in}
\eqdef X_i
\ee
\be
\figplace{dessin3}{0 in}{-0.45 in}
\eqdef
\sum_{i_1,\ldots,i_d}^n
F^{[d]}_{i_1\ldots i_d}
X_{i_1}\ldots X_{i_d}
\ee
Our presentation will now bifurcate in order to specifically address
the three sections of part II.

\subsection{Composition}

Let $F=(F_i)_{1\le i\le n}$
with $F_i\in R[[X_1,\ldots,X_n]]$
be a system of $n$ formal power series without constant
term, given for any $i$, $1\le i\le n$, under the tensorial form
\be
F=\sum_{d\ge 1}\frac{1}{d!}
\sum_{j_1,\ldots,j_d}^n
F^{[d]}_{i,j_1\ldots j_d}
X_{j_1}\ldots X_{j_d}
\ee
We will denote
\be
\figplace{dessin4}{0 in}{-0.45 in}
\eqdef
F^{[d]}_{i,j_1\ldots j_d}
=
\left.
\frac{\partial^d}{\partial X_{j_1}\ldots\partial X_{j_d}}
\right|_{X=0}
F_i(X)
\ee
$F^{[d]}_{i,j_1\ldots j_d}$
can be thought of as a tensor element of a $1$-contravariant
and $d$-covariant tensor, which is symmetric in the $d$ covariant indices
i.e. an element of the $R$-module
$R^n\otimes Sym^d((R^n)^\ast)$.
Suppose we have two systems of $n$ formal power series without constant term,
$F=(F_i)_{1\le i\le n}$ and $G=(G_i)_{1\le i\le n}$.
Our purpose in this section is to give a precise mathematical meaning
to formula (\ref{FrondG}).
For this we need to introduce some definitions and notations.
Let 
${\Br s}=({\Br s}_1,\ldots,{\Br s}_n)$,
$s=(s_1,\ldots,s_n)$,
${\Br t}=({\Br t}_1,\ldots,{\Br t}_n)$,
$t=(t_1,\ldots,t_n)$,
${\Br u}=({\Br u}_1,\ldots,{\Br u}_n)$
and
$u=(u_1,\ldots,u_n)$ be  six vectors of indeterminates.
Denote by $K$ the ring of formal power series
$R[[{\Br s},s,{\Br t},t,{\Br u},u]]$
defined over $R$ using these $6n$ indeterminates.
An element $U\in K$ can be written in multiindex notation as
\be
U=\sum_{\al_1,\ldots,\al_6\in\NN^n}
u_{\al_1,\ldots,\al_6}
\frac{{\Br s}^{\al_1}s^{\al_2}
{\Br t}^{\al_3}t^{\al_4}
{\Br u}^{\al_5}u^{\al_6}}
{\al_1!\ldots\al_6!}
\ee
but also in tensorial notation as
\be
U=\sum_{k_1,\ldots,k_6\in\NN}
\sum_{\ta_1,\ldots,\ta_6}
\frac{U_{\ta_1,\ldots,\ta_6}^{[k_1,\ldots,k_6]}}{k_1!\ldots k_6!}
{\Br s}_{\ta_1}s_{\ta_2}
{\Br t}_{\ta_3}t_{\ta_4}
{\Br u}_{\ta_5}u_{\ta_6}
\ee
where the sum on $\ta_1$ is over all maps $[k_1]\rightarrow [n]$
and likewise for $\ta_2,\ldots,\ta_6$.
${\Br s}_{\ta_1}\eqdef{\Br s}_{\ta_1(1)}\ldots{\Br s}_{\ta_1(k_1)}$
and likewise for $s_{\ta_2}$,
${\Br t}_{\ta_3}$, $t_{\ta_4}$,
${\Br u}_{\ta_5}$ and $u_{\ta_6}$.
Given $A$, $B$ and $C$ three matrices in $GL_n(R)$
we define
\bea
\lefteqn{
I_{A,B,C}(\ta_1,\ldots,\ta_6)\eqdef
\sum_{\si,\mu,\nu}
\lp
\prod_{1\le k\le k_1}
[A^{-1}]_{\ta_2(\si(k))\ta_1(k)}
\rp
} & & \nonumber \\
 & & \times
\lp
\prod_{1\le k\le k_3}
[B^{-1}]_{\ta_4(\mu(k))\ta_3(k)}
\rp
\lp
\prod_{1\le k\le k_5}
[C^{-1}]_{\ta_6(\nu(k))\ta_5(k)}
\rp
\eea
where the sum is over all bijective maps $\si:[k_1]\rightarrow[k_2]$,
$\mu:[k_3]\rightarrow[k_4]$ and
$\nu:[k_5]\rightarrow[k_6]$. Of course, the result is zero
unless $k_1=k_2$, $k_3=k_4$ and $k_5=k_6$.

\begin{lemma}
Let $\rh_1,\ldots,\rh_6$ be in $\GS_{k_1},\ldots,\GS_{k_6}$
respectively, then
\be
I_{A,B,C}(\ta_1\circ\rh_1,\ldots,\ta_6\circ\rh_6)
=
I_{A,B,C}(\ta_1,\ldots,\ta_6)
\ee
\end{lemma}

\noindent{\bf Proof :}
Trivial.
\endproof

\begin{truc}
Let $\al_1,\ldots,\al_6\in\NN^n$,
we define the {\em formal Gaussian integral} of the monomial
${\Br s}^{\al_1}s^{\al_2} {\Br t}^{\al_3}t^{\al_4}
{\Br u}^{\al_5}u^{\al_6}$
, with covariances $A^{-1}$, $B^{-1}$, $C^{-1}$, as
\[
(det\ A)^{-1}
(det\ B)^{-1}
(det\ C)^{-1}
I_{A,B,C}(\ta_1,\ldots,\ta_6)
\]
which belongs to $R$, and where for each $i$, $1\le i\le 6$,
$\ta_i$ is any map from $[k_i]$ to $[n]$ with $k_i=|\al_i|$
and such that the associated multiplicity multiindex
$\mu(\ta_i)=\al_i$.
By the previous lemma, this element of $R$ is independent
of the choice of $\ta_1,\ldots,\ta_6$.
We {\em denote} this expression by
\[
\int d{\Br s}ds
d{\Br t}dt d{\Br u}du
\ e^{-{\Br s}As-{\Br t}Bt-{\Br u}Cu}
\ {\Br s}^{\al_1}s^{\al_2} {\Br t}^{\al_3}t^{\al_4}
{\Br u}^{\al_5}u^{\al_6}
\]
\end{truc}

\begin{truc}
If $U$ is as before a power series in $K$, we define the {\em formal Gaussian
integral}
of $U$ with covariances $A^{-1}$, $B^{-1}$ and $C^{-1}$ as
\bea
\lefteqn{
\int d{\Br s}ds
d{\Br t}dt d{\Br u}du
\ e^{-{\Br s}As-{\Br t}Bt-{\Br u}Cu}
\ U
\eqdef
\sum_{\al_1,\ldots,\al_6\in\NN^n}
\frac{u_{\al_1,\ldots,\al_6}}
{\al_1!\ldots\al_6!}
} & & \nonumber\\
& &
\int d{\Br s}ds
d{\Br t}dt d{\Br u}du
\ e^{-{\Br s}As-{\Br t}Bt-{\Br u}Cu}
\ {\Br s}^{\al_1}s^{\al_2} {\Br t}^{\al_3}t^{\al_4}
{\Br u}^{\al_5}u^{\al_6}
\eea
if the right hand side is summable.
$R$ being equipped with the discrete topology, this simply means that
finitely many terms are nonzero.
\end{truc}

We now restrict to the case where $A=B=C=I$ the $n\times n$
identity matrix, and can state
\begin{theor}
For any $d\ge 1$ and $i,j_1,\ldots,j_d\in [n]$,
\be
(F\circ G)^{[d]}_{i,j_1\ldots j_d}=
\int d{\Br s}ds
d{\Br t}dt d{\Br u}du
\ e^{-{\Br s}s-{\Br t}t-{\Br u}u}
\ U
\ee
where
\be
U\eqdef
s_i{\Br u}_{j_1}\ldots{\Br u}_{j_d}
\exp\lp
{\Br s}F(t)+{\Br t}G(u)
\rp
\ee
\end{theor}

\medskip
\noindent{\bf Remark :}
Again we used our tensorial notation for the coefficients of
$F\circ G(X)$, that is
\be
(F\circ G)^{[d]}_{i,j_1\ldots j_d}
\eqdef
\left.
\frac{\partial^d}{\partial X_{j_1}\ldots\partial X_{j_d}}
\right|_{X=0}
F_i(G(X))
\ee

\medskip
\noindent{\bf Remark :}
As ${\Br s}F(t)$ and ${\Br t}G(u)$
have no constant term, $U$ is a well-defined
element of $K$.

\medskip
\noindent{\bf Remark :}
If one formally expands (\ref{FrondG}) with respect to $X$ one obtains
\bea
\lefteqn{
\sum_{d\ge 1}\sum_{j_1,\ldots,j_d}^n
\frac{X_{j_1}\ldots X_{j_d}}{d!}
(F\circ G)^{[d]}_{i,j_1\ldots j_d}
=\sum_{d\ge 1}\sum_{j_1,\ldots,j_d}^n
\frac{X_{j_1}\ldots X_{j_d}}{d!}
} & & \nonumber \\
 & &
\int d{\Br s}ds
d{\Br t}dt d{\Br u}du
\ e^{-{\Br s}s-{\Br t}t-{\Br u}u}
s_i{\Br u}_{j_1}\ldots{\Br u}_{j_d}
e^{{\Br s}F(t)+{\Br t}G(u)}
\eea
Therefore Theorem 2 is a rigorous restatement
of Claim 1.

\medskip
In order to prove the theorem, we need to define the notions
of pre-Feynman and Feynman diagram structures, which find their
natural habitat in the Joyal theory
of combinatorial species.
The proof of Theorem 2 will accordingly be postponed till the
end of this section.
We suppose that 
$d\ge 1$. The index $i$ considered as a map from $I\eqdef [1]$
to $[n]$ and the collection of indices
$j_1,\ldots,j_d$ considered as a map from $J\eqdef [d]$
to $[n]$ are {\em fixed} in the following.
These two maps we call {\em index assignments}.
Now let $E$ be any finite set.
\begin{truc}
A pre-Feynman diagram sructure on $E$
is an ordered collection
\be
\cE=(E_{\Br s}, E_s, E_{\Br t}, E_t, E_{\Br u}, E_u,
E_{int}, E_{ext}, \pi_F,\pi_G, \rh_s,\rh_{\Br u})
\ee
made of the following data.
\begin{itemize}
\item
$E_{\Br s}, E_s, E_{\Br t}, E_t, E_{\Br u}, E_u,
E_{int}, E_{ext}$
are subsets of $E$.
\item
$\pi_F,\pi_G$ are (unordered) sets of subsets of $E$.
\item
$\rh_s$ is a map from $I$ to $E_{ext}\cap E_s$.
\item
$\rh_{\Br u}$ is a map from $J$ to $E_{ext}\cap E_{\Br u}$.
\end{itemize}

We furthermore ask that the previous data satisfy the
following constraints.
\begin{itemize}
\item
$E_{\Br s}, E_s, E_{\Br t}, E_t, E_{\Br u}, E_u$ are disjoint and their
union is $E$.
\item
$E_{int}, E_{ext}$ are disjoint and their
union is $E$.
\item
$(E_{\Br s}\cup E_{\Br t}\cup E_t\cup E_u)\cap E_{ext}=\emptyset$
\item
$\rh_s:I\rightarrow E_{ext}\cap E_s$
and
$\rh_{\Br u}:J\rightarrow E_{ext}\cap E_{\Br u}$
are bijective.
\item
$\pi_F\cap\pi_G=\emptyset$ and
$\pi_F\cup\pi_G$ forms a partition of $E_{int}$.
\item
Any block $B\in\pi_F$, also called an $F$-vertex,
is the union of $B\cap E_{\Br s}$ and $B\cap E_{t}$
which must respectively be a singleton and a nonempty set.
\item
Any block $B\in\pi_G$, also called a $G$-vertex,
is the union of $B\cap E_{\Br t}$ and $B\cap E_{u}$
which must respectively be a singleton and a nonempty set.
\end{itemize}
\end{truc}

We denote the set of pre-Feynman diagram structures on a finite set $E$ by
$PreFey(E)$ which is obviously finite too.
What we have just done is defining a covariant
endofunctor for the groupoid category of finite sets with morphisms
given by bijective maps.
Indeed, if $\si:E\rightarrow E'$ is such a morphism, its transform by
this functor $PreFey(\si):PreFey(E)\rightarrow PreFey(E')$
is the map which to a pre-Feynman diagram structure
\be
\cE=(E_{\Br s}, E_s, E_{\Br t}, E_t, E_{\Br u}, E_u,
E_{int}, E_{ext}, \pi_F,\pi_G, \rh_s,\rh_{\Br u})
\ee
on $E$ associates the analogous structure
\be
\cE'=(E'_{\Br s}, E'_s, E'_{\Br t}, E'_t, E'_{\Br u}, E'_u,
E'_{int}, E'_{ext}, \pi'_F,\pi'_G, \rh'_s,\rh'_{\Br u})
\ee
on $E'$ given in the obvious manner by
$E'_{\Br s}=\si(E_{\Br s})$,
$E'_{s}=\si(E_{s})$,
$E'_{\Br t}=\si(E_{\Br t})$,
$E'_{t}=\si(E_{t})$,
$E'_{\Br u}=\si(E_{\Br u})$,
$E'_{u}=\si(E_{u})$,
$E'_{int}=\si(E_{int})$,
$E'_{ext}=\si(E_{ext})$,
$\pi'_F=\{\si(B)|B\in\pi_F\}$,
$\pi'_G=\{\si(B)|B\in\pi_G\}$,
$\rh'_s=\si\circ\rh_s$ and finally
$\rh'_{\Br u}=\si\circ\rh_{\Br u}$.
We have thus constructed an example of combinatorial species
in the sense of Joyal, denoted by $PreFey$.
$PreFey(\si)$ is the transport of structure along the bijection $\si$
between the finite sets $E$ and $E'$.

As the previous definition might be hard to digest if served dry,
let us pause to explain the rationale and give
an example.
The ``job'' of a pre-Feynman diagram structure $\cE$ on a set $E$
is to encode an algebraic formula.
We have so to speak defined a ``programming language'' with 
its syntactic rules (the constraints in Definition 3);
a ``program'' in this language (a pre-Feynman diagram structure)
serves to compute an element of the ring $K$.
For example take $d=5$ and $E=[16]$,
with
$E_{\Br s}=\{2\}$,
$E_{s}=\{1\}$,
$E_{\Br t}=\{5,9\}$,
$E_{t}=\{3,4\}$,
$E_{\Br u}=\{12,13,14,15,16\}$,
$E_{u}=\{6,7,8,10,11\}$,
$E_{int}=\{2,3,4,5,6,7,8,9,10,11\}$,
$E_{ext}=\{1,12,13,14,15,16\}$,
$\pi_F=\{\{2,3,4\}\}$,
$\pi_G=\{\{5,6,7,8\}\{9,10,11\}\}$,
$\rh_s:[1]\rightarrow E_s\cap E_{ext}=\{1\}$
(two a priori unrelated copies of the set $\{1\}$)
is given by $\rh_s(1)=1$, and finally let
$\rh_{\Br u}:[d]\rightarrow E_{\Br u}\cap E_{ext}=
\{12,13,14,15,16\}$ be given by
$\rh_{\Br u}(1)=12$,
$\rh_{\Br u}(2)=13$,
$\rh_{\Br u}(3)=14$,
$\rh_{\Br u}(4)=15$ and
$\rh_{\Br u}(5)=16$.
The element of $K$ computed by this data
or ``program'' is
\[
s_i
{\Br u}_{j_1}
{\Br u}_{j_2}
{\Br u}_{j_3}
{\Br u}_{j_4}
{\Br u}_{j_5}
\times
\lp
\sum_{\al_2,\al_3,\al_4=1}^n
{\Br s}_{\al_2}
F^{[2]}_{\al_2,\al_3\al_4}
t_{\al_3}
t_{\al_4}
\rp
\]
\[
\times
\lp
\sum_{\al_5,\al_6,\al_7,\al_8=1}^n
{\Br t}_{\al_5}
G^{[3]}_{\al_5,\al_6\al_7\al_8}
u_{\al_6}
u_{\al_7}
u_{\al_8}
\rp
\]
\be
\times
\lp
\sum_{\al_9,\al_{10},\al_{11}=1}^n
{\Br t}_{\al_9}
G^{[2]}_{\al_9,\al_{10}\al_{11}}
u_{\al_{10}}
u_{\al_{11}}
\rp
\label{expreF}
\ee
A quick look a this expression will convince the reader that a much better
way to represent it is by the following picture
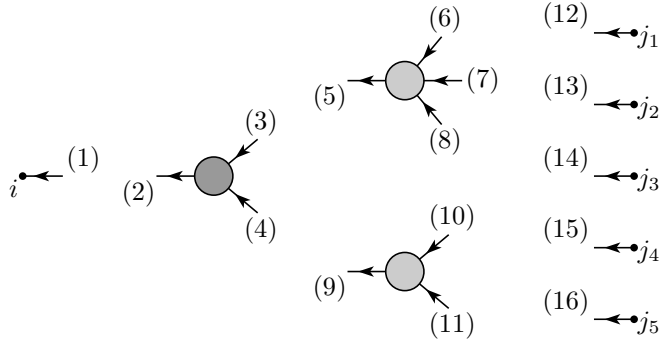
\[
\figput{dessin5}
\]
Remember that we have fixed in our discussion $d$, $i$ and
$j_1,\ldots,j_5$.
This piece of information that is necessary to write and
evaluate the expression (\ref{expreF}) is not included
in the data carried by our pre-Feynman diagram
structure.
The set $E$, in our example serves as an {\em abstract} set
of labels for the indeterminates of type
$\Br s$, $s$,
$\Br t$, $t$,
$\Br u$ and  $u$ that appear in (\ref{expreF}) and are represented as
oriented half-lines in the picture where we have
indicated the labelling between parentheses.
These indeterminates are called ``fields'' in the QFT terminology.
The expression (\ref{expreF}) is called the {\em amplitude}
of the previous pre-Feynman diagram structure given the
additional external structure $i$, $j_1,\ldots,j_5\in [n]$.
We now can introduce the notion of a Feynman diagram structure.
\begin{truc}
A Feynman diagram structure on a finite set $E$ is a quadruple
$\cF=(\cE,\cC_s,\cC_t,\cC_u)$
made of a pre-Feynman diagram structure
\be
\cE=(E_{\Br s}, E_s, E_{\Br t}, E_t, E_{\Br u}, E_u,
E_{int}, E_{ext}, \pi_F,\pi_G, \rh_s,\rh_{\Br u})
\ee
and three {\em bijective} maps
$\cC_s:E_{\Br s}\rightarrow E_s$,
$\cC_t:E_{\Br t}\rightarrow E_t$ and
$\cC_u:E_{\Br u}\rightarrow E_u$.

The maps $\cC_s$, $\cC_t$, $\cC_u$ are called {\em contraction schemes}
for the
$s$-fields, $t$-fields and $u$-fields respectively.
\end{truc}
The set of Feynman diagram structures on $E$ is denoted
by $Fey(E)$.
We are again defining a functor which is a combinatorial specie.
Indeed, if $\si:E\rightarrow E'$ is a bijective map, we define
the transformed morphism $Fey(\si):Fey(E)\rightarrow Fey(E')$
in the obvious manner by letting
$Fey(\si)(\cE,\cC_s,\cC_t,\cC_u)=(\cE',\cC'_s,\cC'_t,\cC'_u)$
with $\cE'=PreFey(\cE)$,
$\cC'_s=\si\circ\cC_s\circ(\si^{-1})|_{E'_{\Br s}}$,
$\cC'_t=\si\circ\cC_t\circ(\si^{-1})|_{E'_{\Br t}}$ and
$\cC'_u=\si\circ\cC_u\circ(\si^{-1})|_{E'_{\Br u}}$.

Again the idea is to encode, thanks to such a structure, an algebraic
expression whose value lies this time in the ground ring $R$ instead
of the formal power series ring $K$.
For instance, if we take the previous example of  pre-Feynman
diagram structure and add to it the maps
$\cC_s$, $\cC_t$, $\cC_u$ given by
$\cC_s(2)=1$,
$\cC_t(5)=3$,
$\cC_t(9)=4$,
$\cC_u(12)=6$,
$\cC_u(13)=7$,
$\cC_u(14)=8$,
$\cC_u(15)=10$ and
$\cC_u(16)=11$;
the resulting algebraic expression,
also called the amplitude of this Feynman diagram structure, is
\[
\sum_{\al_2,\ldots,\al_{11}}^n
\de_{i\al_2}
F^{[2]}_{\al_2,\al_3\al_4}
\de_{\al_3\al_5}
G^{[3]}_{\al_5,\al_6\al_7\al_8}
\]
\[
\de_{\al_6 j_1}
\de_{\al_7 j_2}
\de_{\al_8 j_3}
\de_{\al_4\al_9}
G^{[2]}_{\al_9,\al_{10}\al_{11}}
\de_{\al_{10} j_4}
\de_{\al_{11} j_5}
\]
which belongs to $R$ and where $\de_{ij}$
is simply Kronecker's symbol whose presence is due to our choice
$A=B=C=I$ defining the covariance matrices.
Again it does not take long to realize that a much better notation
for this messy formula
is
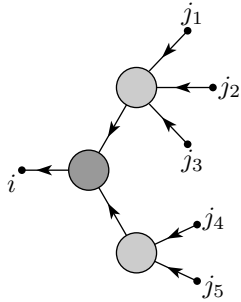
\[
\figput{dessin6}
\]
Note that this expression becomes, when varying the indices $i$,
$j_1,\ldots,j_5$ in $[n]$, the collection of entries of a tensor
in $R\otimes(R^\ast)^{\otimes 5}$
that is built using tensor contraction from three elementary tensors
corresponding to some homogenous components of the series
$F$ and $G$.
We have thus combinatorially translated a very natural
``conceptual'' construction of multilinear algebra.
Note that this composite tensor has no reason to be
symmetric in
$j_1,\ldots,j_5$.

Having provided the definitions of pre-Feynman and Feynman
diagram structures, and an example illustrating their meaning,
we will now, for the sake of mathematical precision,
give the formal definition of amplitudes.

\begin{truc}
Let as before $\cE$ be a pre-Feynman diagram
structure on a finite set $E$, and suppose we are
given two
assignment maps $\ta_s:I\rightarrow [n]$
and $\ta_{\Br u}:J\rightarrow [n]$
with $I=[1]$, $J=[d]$.
We call an {\em index attribution}  any map
$\al:E\rightarrow [n]$ such that $\al|_{E_{ext}\cap E_s}=
\ta_s\circ\rh_s^{-1}$ and $\al|_{E_{ext}\cap E_{\Br u}}=
\ta_{\Br u}\circ\rh_{\Br u}^{-1}$.
Given such an index attribution map $\al$ and a block
$B\in\pi_F$,
if $B\cap E_{\Br s}=\{x\}$
and
$B\cap E_{t}=\{y_1,\ldots,y_p\}$
with $p\ge 1$ we denote
\be
F(B,\al)\eqdef
F^{[p]}_{\al(x),\al(y_1)\ldots\al(y_p)}
\ee
which does not depend on the chosen order of the elements
in $B\cap E_{t}$.
Likewise,
if
$B\in\pi_G$ is such that
$B\cap E_{\Br t}=\{x\}$
and
$B\cap E_{u}=\{y_1,\ldots,y_p\}$
with $p\ge 1$ we denote
\be
G(B,\al)\eqdef
G^{[p]}_{\al(x),\al(y_1)\ldots\al(y_p)}
\ee
We now define the {\em amplitude} of the pre-Feynman
diagram structure $\cE$ on $E$ with respect to the
assignment
maps $\ta_s$ and $\ta_{\Br u}$
as
\bea
\lefteqn{
\cA_{PreFey}(E,\cE,\ta_s,\ta_{\Br u})\eqdef
\sum_{\al}
\lp
\prod_{{\Br x}\in E_{\Br s}}{\Br s}_{\al({\Br x})}
\rp
\lp
\prod_{{x}\in E_{s}}{s}_{\al({x})}
\rp
} & & \nonumber\\
 & &
\lp
\prod_{{\Br y}\in E_{\Br t}}{\Br t}_{\al({\Br y})}
\rp
\lp
\prod_{{y}\in E_{t}}{t}_{\al({y})}
\rp
\lp
\prod_{{\Br z}\in E_{\Br u}}{\Br u}_{\al({\Br z})}
\rp
\lp
\prod_{{z}\in E_{u}}{u}_{\al({z})}
\rp
\nonumber\\
 & & 
\lp
\prod_{B\in\pi_F} F(B,\al)
\rp
\lp
\prod_{B\in\pi_G} G(B,\al)
\rp
\eea
where the sum is over all index attribution maps $\al$.
$\cA_{PreFey}(E,\cE,\ta_s,\ta_{\Br u})$ belongs to the ring
$K=R[[{\Br s},s,{\Br t},t,{\Br u},u]]$.
\end{truc}

\begin{truc}
With the same notation as before, to a Feynman diagram structure
$\cF$ on $E$ and two assignment maps $\ta_s$ and $\ta_{\Br u}$
we associate the corresponding {\em amplitude}
\bea
\lefteqn{
\cA_{Fey}(E,\cF,\ta_s,\ta_{\Br u})
\eqdef
\sum_\al
\lp
\prod_{{\Br x}\in E_{\Br s}}
\de_{\al(\cC_s(\Br x))\al({\Br x})}
\rp
\lp
\prod_{{\Br y}\in E_{\Br t}}
\de_{\al(\cC_t(\Br y))\al({\Br y})}
\rp
} & & \nonumber \\
 & & \times
\lp
\prod_{{\Br z}\in E_{\Br u}}
\de_{\al(\cC_u(\Br z))\al({\Br z})}
\rp
\lp
\prod_{B\in\pi_F} F(B,\al)
\rp
\lp
\prod_{B\in\pi_G} G(B,\al)
\rp
\eea
where again the sum is over all index attribution maps $\al$
compatible with the external structure provided by
$\ta_s$ and $\ta_{\Br u}$, and $\de_{ij}$
is the Kronecker symbol.
Note that this time $\cA_{Fey}(E,\cF,\ta_s,\ta_{\Br u})$
belongs to $R$.
\end{truc}

The following important propositions are obvious from the previous
definitions, and state the relabelling invariance of the amplitudes.

\begin{prop}
If $E$, $E'$ are two finite sets equipped with pre-Feynman diagram
structures $\cE$ and $\cE'$ respectively, such that there exists
a bijection $\si:E\rightarrow E'$ that sends $\cE$ on $\cE'$
by $PreFey(\si)$,
then
\be
\cA_{PreFey}(E,\cE,\ta_s,\ta_{\Br u})=
\cA_{PreFey}(E',\cE',\ta_s,\ta_{\Br u})
\ee
\end{prop}

\begin{prop}
If $E$, $E'$ are two finite sets equipped with Feynman diagram
structures $\cF$ and $\cF'$ respectively, such that there exists
a bijection $\si:E\rightarrow E'$ that sends $\cF$ on $\cF'$
by $Fey(\si)$,
then
\be
\cA_{Fey}(E,\cF,\ta_s,\ta_{\Br u})=
\cA_{Fey}(E',\cF',\ta_s,\ta_{\Br u})
\ee
\end{prop}

An important notion is that of {\em automorphism group}
(the group of {\em ambiguity} we mentioned earlier)
of pre-Feynman and Feynman diagram structures.

\begin{truc}
The automorphism group of a pre-Feynman diagram structure
$\cE$ on $E$ is the group $Aut(E,\cE)$ of all bijective maps
$\si:E\rightarrow E'$ such that
$PreFey(\si)$ leaves $\cE$ unchanged.
\end{truc}

\begin{prop}
\bea
\lefteqn{
\#(Aut(E,\cE))=
} & & \nonumber \\
 & & 
\prod_{p\ge 1}
\lp
m_{F,p}! (p!)^{m_{F,p}}
\rp
\times
\prod_{q\ge 1}
\lp
m_{G,q}! (q!)^{m_{G,q}}
\rp
\eea
where for each integer $p\ge 1$,
$m_{F,p}$ is the number of blocks $B\in\pi_F$ such that
$\#(B\cap E_t)=p$
and for each  integer $q\ge 1$,
$m_{G,q}$ is the number of blocks $B\in\pi_G$ such that
$\#(B\cap E_u)=q$.
\end{prop}

\noindent{\bf Proof :}
A map $\si:E\rightarrow E$ that preserves the structure $\cE$
is necessarily the identity on
$E_{ext}=(E_{ext}\cap E_s)\cup(E_{ext}\cap E_{\Br u})$,
since the injective maps $\rh_s$ and $\rh_{\Br u}$ satisfy
$\rh_s\circ\si=\rh_s$ on $I$
and $\rh_{\Br u}\circ\si=\rh_{\Br u}$ on $J$.
Besides $\si$ must permute
the blocks of $\pi_F$
that contain the same number of
elements from $E_t$, which accounts for the $m_{F,p}!$ factors.
Likewise $\si$ must permute
the blocks of $\pi_G$
that contain the same number of
elements from $E_u$, which accounts for the $m_{G,q}!$ factors.
Finally $\si$ permutes the elements within each block of $\pi_F$ and
$\pi_G$, which gives the $p!$ and $q!$ factors.
\endproof

\begin{prop}
With the notations of Theorem 2,
let $i$, $j_1,\ldots,j_d$ be elements of $[n]$
that define assignment maps
$\ta_s:I=[1]\rightarrow [n]$
by $\ta_s(1)=i$
and $\ta_{\Br u}:J=[d]\rightarrow [n]$
by $\ta_{\Br u}(\nu)=j_\nu$, for $1\le \nu\le d$.
The element
\[
U=
s_i{\Br u}_{j_1}\ldots{\Br u}_{j_d}
\exp\lp
{\Br s}F(t)+{\Br t}G(u)
\rp
\]
of the ring $K$ can be rewritten as
\be
U=\sum_{[E,\cE]}
\frac{\cA_{preFey}(E,\cE,\ta_s,\ta_{\Br u})}
{\#Aut(E,\cE)}
\label{UpreF}
\ee
where the sum is over the isomorphism classes of pairs $(E,\cE)$
made of a finite set $E$ and a pre-Feynman diagram
structure $\cE$ on $E$.
The sets $I=[1]$, $J=[d]$
and the assignment maps
$\ta_s:I\rightarrow [n]$ and
$\ta_{\Br u}:J\rightarrow [n]$ 
are fixed throughout.
Two pairs $(E,\cE)$ and $(E',\cE')$ are said isomorphic
if there exists a bijection $\si:E\rightarrow E'$
such that $PreFey(\si)$ sends $\cE$
to $\cE'$.
In the sum (\ref{UpreF}) $(E,\cE)$ denotes any representative
of the class $[E,\cE]$.
\end{prop}

\noindent{\bf Proof :}
If $m\ge 1$, we use the shorthand notation
\be
{\Br s}Ft^m\eqdef
\sum_{\al,\beta_1,\ldots,\beta_m=1}^n
{\Br s}_\al
F^{[m]}_{\al,\beta_1\ldots\beta_m}
t_{\beta_1}\ldots t_{\beta_m}
\ee
and
\be
{\Br t}Gu^m\eqdef
\sum_{\al,\beta_1,\ldots,\beta_m=1}^n
{\Br t}_\al
G^{[m]}_{\al,\beta_1\ldots\beta_m}
u_{\beta_1}\ldots u_{\beta_m}
\ee
Therefore in the formal power series ring $K$
we have
\bea
U & = & s_i {\Br u}_{j_1}\ldots {\Br u}_{j_d}
\exp
\lp
\sum_{p\ge 1}\frac{{\Br s}F t^p}{p!}
+
\sum_{q\ge 1}\frac{{\Br t}G u^q}{q!}
\rp \\
 & = &
s_i {\Br u}_{j_1}\ldots {\Br u}_{j_d}
\sum_{m_F,m_G\ge 0}
\frac{1}{m_F! m_G!}
{\lp
\sum_{p\ge 1}\frac{{\Br s}F t^p}{p!}
\rp}^{m_F}
{\lp
\sum_{q\ge 1}\frac{{\Br t}G u^q}{q!}
\rp}^{m_G}
\eea
which by the multinomial theorem becomes
\bea
\lefteqn{
U=\sum_{m_F,m_G\ge 0}
\sum_{(m_{F,p})_{p\ge 1}}
\sum_{(m_{G,q})_{q\ge 1}}
s_i {\Br u}_{j_1}\ldots {\Br u}_{j_d}
} & & \nonumber \\
& & 
\times
\prod_{p\ge 1}
\lp
\frac{1}{m_{F,p}!}
{\lp
\frac{{\Br s}F t^p}{p!}
\rp}^{m_{F,p}}
\rp
\times
\prod_{q\ge 1}
\lp
\frac{1}{m_{G,q}!}
{\lp
\frac{{\Br t}G u^q}{q!}
\rp}^{m_{G,q}}
\rp
\eea
where the sum is over all families $(m_{F,p})_{p\ge 1}$
and $(m_{G,q})_{q\ge 1}$
of nonnegative integers, necessarily of finite support,
such that
$\sum_{p\ge 1}m_{F,p}=m_F$ and
$\sum_{q\ge 1}m_{G,q}=m_G$.
One can then ``remove the parentheses'' in the
above packet summation in $K$, since
for each monomial in the variables
${\Br s}$, $s$,
${\Br t}$, $t$,
${\Br u}$ and $u$, only finitely many $m_F$'s and $m_G$'s contribute.
Then
\bea
\lefteqn{
U=
\sum_{(m_{F,p})_{p\ge 1},(m_{G,q})_{q\ge 1}}
s_i {\Br u}_{j_1}\ldots {\Br u}_{j_d}
} & & \nonumber \\
& & 
\times
\prod_{p\ge 1}
\lp
\frac{1}{m_{F,p}!}
{\lp
\frac{{\Br s}F t^p}{p!}
\rp}^{m_{F,p}}
\rp
\times
\prod_{q\ge 1}
\lp
\frac{1}{m_{G,q}!}
{\lp
\frac{{\Br t}G u^q}{q!}
\rp}^{m_{G,q}}
\rp
\eea
where the sum is over all pairs
$((m_{F,p})_{p\ge 1},(m_{G,q})_{q\ge 1})$
of finitely supported families of nonnegative integers.
Now notice that such pairs are in bijective correspondance
with equivalence classes $[E,\cE]$
of finite sets equipped with pre-Feynman diagram structures.
The previous proposition does the rest.
\endproof

\noindent
We now have as an immediate consequence the following

\begin{prop}
Keeping the same notation, we have in the ring $R$

\bea
\lefteqn{
\int d{\Br s}ds
d{\Br t}dt d{\Br u}du
\ e^{-{\Br s}s-{\Br t}t-{\Br u}u}
\ U=} & & \nonumber \\
 & & \sum_{[E,\cE]}
\sum_{(\cC_s,\cC_t,\cC_u)}
\frac{\cA_{Fey}
\lp
E,(\cE,\cC_s,\cC_t,\cC_u),\ta_s,\ta_{\Br u}
\rp}
{\#Aut(E,\cE)}
\label{contsum}
\eea
where again one sums first over equivalence classes of pre-Feynman
diagram structures, $(E,\cE)$
being an arbitrary representative of such a class.
The $\cC_s$, $\cC_t$, $\cC_u$
are summed over contractions of the $s$, $t$ and $u$ fields respectively
within the chosen representative
$(E,\cE)$.
$\cF=(\cE,\cC_s,\cC_t,\cC_u)$
is then a Feynman diagram structure with underlying pre-Feynman structure
$\cE$.
Also the sum on the right hand side is of finite support.
\end{prop}

\noindent{\bf Proof :}
From Definitions 1 and 2  it is clear that term by term,
i.e.
for each class $[E,\cE]$
we have
\bea
\lefteqn{
\int d{\Br s}ds
d{\Br t}dt d{\Br u}du
\ e^{-{\Br s}s-{\Br t}t-{\Br u}u}
\ \cA_{PreFey}
\lp
E,\cE,\ta_s,\ta_{\Br u}
\rp
=
} & & \nonumber \\
& & 
\sum_{(\cC_s,\cC_t,\cC_u)}
\cA_{Fey}
\lp
E,(\cE,\cC_s,\cC_t,\cC_u),\ta_s,\ta_{\Br u}
\rp
\eea
Proposition 4 being proven, all one has to check is that the sum
over $[E,\cE]$
is well-defined in $R$, i.e.
has finite support.
However this is an easy consequence
of the tree structure of our Feynman diagrams, where the root
is the unique element of $E_s$,
the first generation of vertices corresponds
to the blocks $\pi_F$, the second generation corresponds to those of $\pi_G$
and the third generation, that is the set of leaves of the tree,
is $E_{\Br u}$.
An easy counting argument
using the definition of our pre-Feynman and Feynman diagram structures
shows that in order for a triple of bijections $(\cC_s,\cC_t,\cC_u)$
to exist for a given pair $(E,\cE)$,
one must have
$\#(\pi_F)=1$ and $\#(\pi_G)\le d$,
because $G$ has no constant term.
This shows that finitely many classes $[E,\cE]$ contribute.
\endproof
We now introduce as we did for pre-Feynman diagram structures, the following
definition

\begin{truc}
The automorphism group of a Feynman diagram structure
$\cF$ on $E$ is the group $Aut(E,\cF)$ of all bijective maps
$\si:E\rightarrow E'$ such that
$Fey(\si)$ leaves $\cF$ unchanged.
\end{truc}
We now have the following result
\begin{theor}
With the same notation as before
\be
\int d{\Br s}ds
d{\Br t}dt d{\Br u}du
\ e^{-{\Br s}s-{\Br t}t-{\Br u}u}
\ U
= 
\sum_{[E,\cF]}
\frac{\cA_{Fey}
\lp
E,\cF,\ta_s,\ta_{\Br u}
\rp}
{\#Aut(E,\cF)}
\label{Fexp}
\ee
where the sum is over equivalence classes of pairs $(E,\cF)$
made of a finite set $E$ and a Feynman diagram structure $\cF$ on $E$.
Again the sets
$I=[1]$, $J=[d]$
and the index assignment maps $\ta_s:I\rightarrow [n]$
and $\ta_{\Br u}:J\rightarrow [n]$
are fixed.
Two pairs $(E,\cF)$ and $(E',\cF')$
are said equivalent if there exists a bijection $\si:E\rightarrow E'$
such that $Fey(\si)(\cF)=\cF'$.
In the sum on the right hand side
of (\ref{Fexp}), $(E,\cF)$ denotes any representative of the class
$[E,\cF]$.
\end{theor}
\noindent{\bf Proof :}
Starting from the previous proposition,
one notices that any equivalence class of Feynman diagrams occurs in
equation (\ref{contsum}),
where one takes the equivalence class of
$(E,(\cE,\cC_s,\cC_t,\cC_u))$.
One simply has to count how many times a given Feynman diagram class
appears that way.
If $(E,(\cE,\cC_s,\cC_t,\cC_u))$
is equivalent to
$(E',(\cE',\cC'_s,\cC'_t,\cC'_u))$,
that means that there is a bijection
$\si:E\rightarrow E'$
such that $PreFey(\si)(\cE)=\cE'$,
$\cC'_s=\si\circ\cC_s\circ(\si^{-1})|_{E_{\Br s}}$,
$\cC'_t=\si\circ\cC_t\circ(\si^{-1})|_{E_{\Br t}}$ and
$\cC'_u=\si\circ\cC_u\circ(\si^{-1})|_{E_{\Br u}}$.
But then $(E,\cE)$ is equivalent to $(E',\cE')$
as pre-Feynman diagram structures, which forces
$E=E'$ and $\cE=\cE'$,
since in equation (\ref{contsum}) one takes only one representative in
each class of pre-Feynman diagram structures.
Therefore the number of occurences in (\ref{contsum}) of the
class of $(E,\cF)$, with $\cF=(\cE,\cC_s,\cC_t,\cC_u)$,
is equal to the number of contractions
$(\cC'_s,\cC'_t,\cC'_u)$ such that
$(E,(\cE,\cC'_s,\cC'_t,\cC'_u))$ is equivalent to
$(E,(\cE,\cC_s,\cC_t,\cC_u))$.
In other words, one is counting the cardinality of the
orbit of $(\cC_s,\cC_t,\cC_u)$
under the left-action
of $Aut(E,\cE)$
given by
\be
\si.(\cC_s,\cC_t,\cC_u)
=
\lp
\si\circ\cC_s\circ(\si^{-1})|_{E_{\Br s}},
\ \si\circ\cC_t\circ(\si^{-1})|_{E_{\Br t}},
\ \si\circ\cC_u\circ(\si^{-1})|_{E_{\Br u}}
\rp
\ee
It is equal to
\[
\frac{\#Aut(E,\cE)}{\#Aut(E,\cF)}
\]
since $Aut(E,\cF)$
is the isotropy subgroup of $(\cC_s,\cC_t,\cC_u)$.
Now equation (\ref{Fexp}) follows immediately.
\endproof

\noindent{\bf Remark :}
Equation (\ref{Fexp})
shows that
a formal Gaussian integral
is analogous to a Hurewitz or exponential
generating series (see~\cite{Joyal}).
Indeed it is easy to check that
\be
\sum_{[E,\cF]}
\frac{\cA_{Fey}
\lp
E,\cF,\ta_s,\ta_{\Br u}
\rp}
{\#Aut(E,\cF)}
=
\sum_{k\ge 0}
\sum_{\cF\in Fey([k])}
\frac{\cA_{Fey}
\lp
[k],\cF,\ta_s,\ta_{\Br u}
\rp}{k!}
\label{Hurewitz}
\ee
The classical Hurewitz
series for the species of Feynman diagrams corresponds to
replacing
$\cA_{Fey}(E,\cF,\ta_s,\ta_{\Br u})$ by the coarser
``constant over the orbits'' function
$V^{\#(E)}$ for some indeterminate $V$.

We can now finally proceed to

\medskip
\noindent{\bf Proof of Theorem 2 :}
Let $d\ge 1$ be fixed for the moment, and consider
the polynomial in $R[X_1,\ldots,X_n]$
\bea
\lefteqn{
I_d=
\sum_{\ta_{\Br u}}
X_{\ta_{\Br u}(1)}\ldots X_{\ta_{\Br u}(d)}
\int d{\Br s}ds
d{\Br t}dt d{\Br u}du
\ e^{-{\Br s}s-{\Br t}t-{\Br u}u}
} & & \nonumber \\
& & 
s_i{\Br u}_{\ta_{\Br u}(1)}\ldots{\Br u}_{\ta_{\Br u}(d)}
\exp\lp
{\Br s}F(t)+{\Br t}G(u)
\rp
\eea
where the sum is over all maps
$\ta_{\Br u}:J=[d]\rightarrow [n]$.
We have using Proposition 5, all sums being finite here,
\be
I_d=\sum_{\ta_{\Br u}}
X_{\ta_{\Br u}(1)}\ldots X_{\ta_{\Br u}(d)}
\sum_{[E,\cE]}
\sum_{(\cC_s,\cC_t,\cC_u)}
\frac{\cA_{Fey}
\lp
E,(\cE,\cC_s,\cC_t,\cC_u),\ta_s,\ta_{\Br u}
\rp}
{\#Aut(E,\cE)}
\ee
Now notice that from the symmetry properties of tensor elements of
$F$ and $G$, the Definition 6 of amplitudes and the tree-like description
of the relevent Feynman diagrams given in the proof of
Proposition 5, it is easy to see that 
\be
\sum_{\ta_{\Br u}}
\cA_{Fey}
\lp
E,(\cE,\cC_s,\cC_t,\cC_u),\ta_s,\ta_{\Br u}
\rp
X_{\ta_{\Br u}(1)}\ldots X_{\ta_{\Br u}(d)}
\label{Gamp}
\ee
only depends on $[E,\cE]$
that is on $(m_{G,q})_{q\ge 1}$
the notation being the same as in Proposition 3.
We used the fact
that $m_{F,p}$ vanishes for all $p\ge 1$ except for $p=m$ where
$m\eqdef \sum_{q\ge 1} m_{G,q}$.
Besides one has $\sum_{q\ge 1} q.m_{G,q}=d$.
We denote (\ref{Gamp}) by 
$\Om((m_{G,q})_{q\ge 1})$.
By Proposition 3
\be
\#Aut(E,\cE)
=
m!
\times
\prod_{q\ge 1}
\lp
m_{G,q}! (q!)^{m_{G,q}}
\rp
\ee
Now the number of triples
$(\cC_s,\cC_t,\cC_u)$
of contraction schemes is $1!\times m!\times d!$,
therefore
\bea
I_d & = &
\sum_{(m_{G,q})_{q\ge 1}|\sum_{q\ge 1} q m_{G,q}=d}
\frac{1!\times m!\times d!\times \Om((m_{G,q})_{q\ge 1}) }
{m!\times
\prod_{q\ge 1}
\lp
m_{G,q}! (q!)^{m_{G,q}}
\rp }\\
 & = &
\sum_{(m_{G,q})_{q\ge 1}|\sum_{q\ge 1} q m_{G,q}=d}
\frac{d!\ \Om((m_{G,q})_{q\ge 1}) }
{\prod_{q\ge 1}
\lp
m_{G,q}! (q!)^{m_{G,q}}
\rp }
\eea
Now for a given $m\ge 1$ and $\om=(\om_1,\ldots,\om_m)$
with $\om_i\ge 1$ for all $i$, $1\le i\le m$ and
$\om_1+\cdots+\om_m=d$, we let
$\mu(\om)=(m_{G,q})_{q\ge 1}$ where $m_{G,q}$ counts the number
of indices $i$, $1\le i\le m$ with $\om_i=q$.
It is easy to see that
\be
\Om((m_{G,q})_{q\ge 1})
=
\sum_{\al_1,\ldots,\al_m=1}^n
F^{[m]}_{i,\al_1\ldots\al_m}
\lp G_{\al_1} X^{\om_1}
\rp
\ldots
\lp G_{\al_m} X^{\om_m}
\rp
\ee
where we used the shorthand notation
\be
G_i X^{\nu}\eqdef
\sum_{j_1,\ldots,j_\nu=1}^n
F^{[\nu]}_{i,j_1\ldots j_\nu}
X_{j_1}\ldots X_{j_\nu}
\ee
Since
the number of $\om$'s for which $\mu(\om)$ is equal to a given
$(m_{G,q})_{q\ge 1}$
is by the multinomial theorem
$\frac{m!}{\prod_{q\ge 1} m_{G,q}!}$,
one has
\bea
\lefteqn{
I_d  = 
\sum_{(m_{G,q})_{q\ge 1}|\sum_{q\ge 1} q m_{G,q}=d}
\frac{\prod_{q\ge 1} m_{G,q}!}{m!}
} & & \nonumber \\
 & & \times \sum_{\om|\mu(\om)=(m_{G,q})_{q\ge 1}}
\frac{d!\ \Om(\mu(\om)) }
{\prod_{q\ge 1}
\lp
m_{G,q}! (q!)^{m_{G,q}}
\rp }
\eea
or
\bea
\lefteqn{
I_d  = 
d!\sum_{m\ge 1}
\ \sum_{\om| \om_1+\cdots+\om_m=d}
\ \frac{1}{m!}
} & & \nonumber \\
 & & \times \sum_{\al_1,\ldots,\al_m=1}^n
F^{[m]}_{i,\al_1\ldots\al_m}
\lp G_{\al_1} X^{\om_1}
\rp
\ldots
\lp G_{\al_m} X^{\om_m}
\rp
\eea
and finally summing over $d\ge 1$, we have
\be
\sum_{d\ge 1}\frac{I_d}{d!}
=F_i(G(X))
\ee
since
\be
G_i(X)=\sum_{\nu\ge 1}
\frac{1}{\nu !} G_i X^{\nu}
\ee
The only thing that remains to be checked to prove
Theorem 2 is that the right hand side of equation (\ref{FrondG})
is symmetric with respect to the indices $j_1,\ldots,j_d$,
which is obvious from Lemma 1 and Definitions 1 and 2.
\endproof

\subsection{Reversion}
As in the beginning of section II.2 we let
$F=(F_i)_{1\le i\le n}$
be a system of $n$ formal power series without constant
term in $R[[X_1,\ldots,X_n]]$,
given by
\be
F_i(X)
=\sum_{d\ge 1}\frac{1}{d!}
\sum_{j_1,\ldots,j_d=1}^n
F^{[d]}_{i,j_1\ldots j_d}
X_{j_1}\ldots X_{j_d}
\ee
We will in fact separate the {\em linear} part
\be
L_i(X)
\eqdef
\sum_{j=1}^n
F^{[1]}_{i,j}
X_j
\ee
from the {\em nonlinear} part
\be
H_i(X)\eqdef
-\sum_{d\ge 2}\frac{1}{d!}
\sum_{j_1,\ldots,j_d=1}^n
F^{[d]}_{i,j_1\ldots j_d}
X_{j_1}\ldots X_{j_d}
\ee
so that $F_i(X)=L_i(X)-H_i(X)$.
The linear part $L(u)$,
which becomes quadratic
after contraction with ${\Br u}$ in order to form
${\Br u}L(u)=\sum_{i=1}^n {\Br u}_i L_i(u)$
is called the {\em free} or {\em Gaussian} part in the physics literature.
The remaining terms in the exponential in (\ref{testfunc}) that is
${\Br u}H(u)+{\Br u}Y$
form the {\em interaction} part.

Let  $A\in\cM_n(R)$
be the matrix with entries $A_{ij}\eqdef F^{[1]}_{i,j}$ for
$1\le i,j\le n$.
We will suppose in this section that $A\in GL_n(R)$.
This is a necessary and sufficient condition
for $F=(F_i)_{1\le i\le n}$ to be invertible
for composition of multivariable power series.
Our aim here is to give a precise meaning and rigorous justification for
Claim 3, giving a formula for the compositional inverse $F^{-1}$.
We introduce the vectors ${\Br u}=({\Br u}_1,\ldots,{\Br u}_n)$,
$u=(u_1,\ldots,u_n)$ and $Y=(Y_1,\ldots,Y_n)$
of indeterminates, and will work only
with the formal power series rings $R[[Y]]$ and $R[[{\Br u},u,Y]]$.

Let $k_1$, $k_2$ be in $\NN$, and let
$\ta_1:[k_1]\rightarrow [n]$ and $\ta_2:[k_2]\rightarrow [n]$
be two given maps.
We again define
\be
I_A(\ta_1,\ta_2)\eqdef
\sum_{\si}
\prod_{1\le k\le k_1}
[A^{-1}]_{\ta_2(\si(k))\ta_1(k)}
\ee
where the sum is over all bijective maps $\si:[k_1]\rightarrow [k_2]$.

\begin{truc}
Let $\al_1,\al_2\in\NN^n$, we define the {\em formal Gaussian integral}
of the monomial ${\Br u}^{\al_1}u^{\al_2}$,
with covariance $A^{-1}$ as the element in $R$ given by
\be
\int
d{\Br u}du
\ e^{-{\Br u}Au}
\ {\Br u}^{\al_1}u^{\al_2}
\eqdef
(det\ A)^{-1}
I_A(\ta_1,\ta_2)
\ee
where each $\ta_i$, for $i=1,2$, is any map from
$[k_i]$ to $[n]$
with $k_i=|\al_i|$ and such that the associated multiplicity multiindex
$\mu(\ta_i)$
is equal to $\al_i$.
Again this definition is independent of the choice of 
$\ta_1$ and $\ta_2$.
\end{truc}

\begin{truc}
If
\be
U=\sum_{\al_1,\al_2,\al_3\in\NN^n}
\frac{u_{\al_1,\al_2,\al_3}}{\al_1!\al_2!\al_3!}
\ {\Br u}^{\al_1}u^{\al_2} Y^{\al_3}
\ee
is an element of $R[[{\Br u},u,Y]]$,
we define the {\em formal Gaussian integral} of $U$ as the element
in $R[[Y]]$ given by
\be
\int
d{\Br u}du
\ e^{-{\Br u}Au}
\ U\eqdef
\sum_{\al_1,\al_2,\al_3\in\NN^n}
\frac{u_{\al_1,\al_2,\al_3} Y^{\al_3}}{\al_1!\al_2!\al_3!}
\int d{\Br u}du
\ e^{-{\Br u}Au}
\ {\Br u}^{\al_1}u^{\al_2}
\ee
provided that the right hand side is summable in $R[[Y]]$
(i.e. for any $\al_3\in\NN^n$ there are only finitely many
$(\al_1,\al_2)$'s giving a nonzero contribution).
\end{truc}

\noindent
Now the numerator in Claim 3 can be
interpreted as the application of this definition in the case where
\be
U=u_i\exp\lp
{\Br u}H(u)+{\Br u}Y
\rp
\ee
while the denominator corresponds to the case
\be
U=\exp\lp
{\Br u}H(u)+{\Br u}Y
\rp
\ee
All one has to do is to prove the summability in $R[[Y]]$.
For this we again need the
Feynman diagrammatic machinery. Again we need both the notions of
pre-Feynman and Feynman diagrams.
We let $I$ and $J$ be two fixed finite sets.

\begin{truc}
A pre-Feynman diagram sructure of type $(I,J)$ on a finite set $E$
is an ordered collection
\be
\cE=(E_{\Br u}, E_u,
E_{int}, E_{ext}, \pi_H,\pi_Y, \rh_u,\rh_{\Br u})
\ee
made of the following data.
\begin{itemize}
\item
$E_{\Br u}, E_u,
E_{int}, E_{ext}$
are subsets of $E$.
\item
$\pi_H,\pi_Y$ are (unordered) sets of subsets of $E$.
\item
$\rh_u$ is a map from $I$ to $E_{ext}\cap E_u$.
\item
$\rh_{\Br u}$ is a map from $J$ to $E_{ext}\cap E_{\Br u}$.
\end{itemize}

We also ask that the previous data satisfy the
following constraints.
\begin{itemize}
\item
$E$ is the disjoint union of $E_u$ and
$E_{\Br u}$.
\item
$E$ is the disjoint union of 
$E_{int}$ and $ E_{ext}$.
\item
$\rh_u:I\rightarrow E_{ext}\cap E_u$
and
$\rh_{\Br u}:J\rightarrow E_{ext}\cap E_{\Br u}$
are bijective.
\item
$\pi_H\cap\pi_Y=\emptyset$ and
$\pi_H\cup\pi_Y$ forms a partition of $E_{int}$.
\item
For any block $B\in\pi_H$, also called an $H$-vertex,
$\#(B\cap E_{\Br u})=1$ and
$\#(B\cap E_u)\ge 2$
\item
For any block $B\in\pi_Y$, also called a $Y$-vertex or a leaf,
$\#(B\cap E_{\Br u})=1$ and $B\cap E_u=\emptyset$.
\end{itemize}
\end{truc}

\begin{truc}
A Feynman diagram structure of type $(I,J)$
on a finite set $E$ is a couple $(\cE,\cC)$
made of a pre-Feynman diagram structure
\be
\cE=(E_{\Br u}, E_u,
E_{int}, E_{ext}, \pi_H,\pi_Y, \rh_u,\rh_{\Br u})
\ee
of type $(I,J)$ on $E$
and a bijective map $\cC:E_{\Br u}\rightarrow E_u$.
\end{truc}

\noindent
Transport of structure is defined in the same obvious manner
as in section III.1, which again provides us with two functors
$PreFey$ and $Fey$
which are combinatorial species in the sense of Joyal.

\begin{truc}
Let $\cE$ be a pre-Feynman diagram
structure of type $(I,J)$ on a finite set $E$, and suppose we are
given two
assignment maps $\ta_u:I\rightarrow [n]$
and $\ta_{\Br u}:J\rightarrow [n]$.
We call an {\em index attribution}  any map
$\al:E\rightarrow [n]$ such that $\al|_{E_{ext}\cap E_u}=
\ta_u\circ\rh_u^{-1}$ and $\al|_{E_{ext}\cap E_{\Br u}}=
\ta_{\Br u}\circ\rh_{\Br u}^{-1}$.
Given such an index attribution map $\al$ and a block
$B\in\pi_H$,
if $B\cap E_{\Br u}=\{{\Br x}\}$
and
$B\cap E_{u}=\{y_1,\ldots,y_p\}$
with $p\ge 2$ we denote
\be
H(B,\al)\eqdef
H^{[p]}_{\al({\Br x}),\al(y_1)\ldots\al(y_p)}
=
-F^{[p]}_{\al({\Br x}),\al(y_1)\ldots\al(y_p)}
\ee
which does not depend on the above enumeration of the elements
of $B\cap E_{u}$.
Likewise,
if
$B\in\pi_Y$ is such that
$B=\{{\Br x}\}$, with ${\Br x}\in E_{\Br u}$,
we denote
\be
Y(B,\al)\eqdef
Y_{\al({\Br x})}
\ee
We can now define the {\em amplitude} of the pre-Feynman
diagram structure $\cE$ on $E$ with respect to the
assignment
maps $\ta_u$ and $\ta_{\Br u}$
as
\bea
\lefteqn{
\cA_{PreFey}(E,\cE,\ta_u,\ta_{\Br u})\eqdef
\sum_{\al}
\lp
\prod_{{\Br x}\in E_{\Br u}}{\Br u}_{\al({\Br x})}
\rp
\lp
\prod_{{x}\in E_{u}}{u}_{\al({x})}
\rp
} & & \nonumber\\
 & & \times \lp
\prod_{B\in\pi_H} H(B,\al)
\rp
\lp
\prod_{B\in\pi_Y} Y(B,\al)
\rp
\eea
which belongs to $R[[{\Br u},u,Y]]$.
Again the sum is over all index attribution maps $\al$.
\end{truc}

\begin{truc}
With the same notation as in the previous definition,
to a Feynman diagram structure
$\cF$ of type $(I,J)$
on $E$ and two assignment maps $\ta_u$ and $\ta_{\Br u}$
we associate the corresponding {\em amplitude}
\bea
\lefteqn{
\cA_{Fey}(E,\cF,\ta_u,\ta_{\Br u})
\eqdef
\sum_\al
\lp
\prod_{{\Br x}\in E_{\Br u}}
{(A^{-1})}_{\al(\cC_s(\Br x))\al({\Br x})}
\rp
} & & \nonumber\\
 & & \times
\lp
\prod_{B\in\pi_H} H(B,\al)
\rp
\lp
\prod_{B\in\pi_Y} Y(B,\al)
\rp
\eea
where
$(A^{-1})_{ij}$ denotes the entries of the covariance
matrix $A^{-1}\in GL_n(R)$.
The amplitude
$\cA_{Fey}(E,\cF,\ta_u,\ta_{\Br u})$
is an element in $R[[Y]]$.
\end{truc}

\noindent
Again these amplitudes are obviously invariant by relabelling or
transport of structure.
One defines as in section III.1 the notions of automorphism groups
of pairs $(E,\cE)$
and $(E,\cF)$
with $\cE$ a pre-Feynman diagram structure
and $\cF$ a Feynman diagram structure on $E$.
The following proposition is proved like its sibling from section III.1.

\begin{prop}
If $\cE$ is pre-Feynman diagram structure on $E$,
\be
\#Aut(E,\cE)=
\prod_{p\ge 1}\lp
m_{H,p}!(p!)^{m_{H,p}}
\rp
\ \times\ m_Y!
\ee
where for each $p\ge 2$, $m_{H,p}$
counts the blocks $B\in\pi_H$ such that
$\#(B\cap E_u)=p$
and $m_Y=\#(\pi_Y)$.
\end{prop}

\noindent
One also has by the same arguments as in Proposition 4

\begin{prop}
Given two finite sets $I$ and $J$ and two
index assignment maps $\ta_u$ and $\ta_{\Br u}$
one has in the ring
$R[[{\Br u},u,Y]]$
\be
\lp
\prod_{i\in I}
u_{\ta_u(i)}
\rp
\lp
\prod_{j\in J}
{\Br u}_{\ta_{\Br u}(j)}
\rp
\exp\lp
{\Br u}H(u)+{\Br u}Y
\rp
=
\sum_{[E,\cE]}
\frac{\cA_{preFey}(E,\cE,\ta_u,\ta_{\Br u})}
{\#Aut(E,\cE)}
\ee
\end{prop}

\medskip
Before we state the analog of Proposition 5
and to take care of issues of summability we
have to analyse more closely the Feynman diagram structure
appearing here.
Given such a structure $\cF$ of type $(I,J)$ on $E$,
we can associate to it
an ordinary digraph $G$ on the set ${\til E}$ defined as
the disjoint union of
${\til E}_u$ the set of one-element subsets of $E_{ext}\cap E_u$,
${\til E}_{\Br u}$ the set of one-element subsets
of $E_{ext}\cap E_{\Br u}$,
${\til E}_H=\pi_H$ and
${\til E}_Y=\pi_Y$.
Therefore ${\til E}$ is a partition of $E$.
Now $G$ is the set of ordered pairs $(a,b)$,
with $a,b\in {\til E}$, such that there exist
$x\in a\cap E_{\Br u}$ and
$y\in b\cap E_u$
such that $y=\cC(x)$.
If the link $(a,b)$ is in $G$ we call $a$ its origin and $b$ its
end.
For example for the Feynman diagram represented by the following
picture
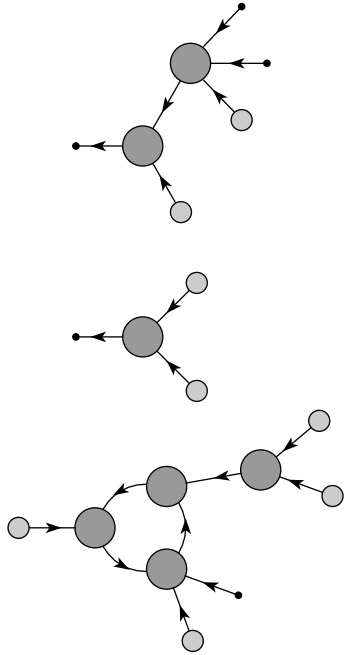
\be
\figput{dessin7}
\ee
$I$ corresponds to the $2$ half-lines
\[
\figput{dessin8}
\]
called the $u$-sources;
$J$ corresponds to the $3$ half-lines
\[
\figput{dessin9}
\]
called the ${\Br u}$-sources.
$E$ is the set of all half-lines and has $2\times 18=36$ elements.
We also have $\#({\til E}_H)=7$,
$\#({\til E}_Y)=8$,
$\#({\til E})=2+7+8+3=20$, and
$\#(G)=18$.
It is a simple but tedious matter of going through the previous
definitions to verify that the only possible connected components of
the digraph $G$ on ${\til E}$ are of two types.

\noindent{\bf Tree-like :}
A tree where all the links are oriented towards the root
that has to be the unique element of ${\til E}_u$ in the component.
The leaves are either $Y$-vertices,
(elements of ${\til E}_Y$) or ${\Br u}$-sources (elements of
${\til E}_{\Br u}$).
The internal vertices of the tree are all $H$-vertices,
i.e. elements of ${\til E}_H$, and {\em have at least two} offsprings.
This crucial property is because $H$
has been defined as the {\em nonlinear}
part of $-F$.

\noindent{\bf Circuit-like :}
A graph with a unique central oriented circuit of $H$-vertices on
which trees like above are hooked. The latter are
oriented towards the circuit,
and their leaves are either $Y$-vertices of ${\Br u}$-sources.
Such a graph contains no element of ${\til E}_u$.

\noindent{\bf Remark :}
Note the analogy with the combinatorial species of endofunctions,
which live here on the ``functorially'' derived
{\em abstract} set ${\til E}$.
No reference is made to the {\em concrete} set of indices $[n]$,
or to the
dimensionality $n$ of the problem, which only appear
in the calculation of amplitudes.
The need of varying $n$, in order to realize
the manifold $\cB$ mentioned in the introduction as an
``inductive limit'' of finite sets (and thus the set of maps
$\cB\rightarrow\cT$ as a ``projective limit''), makes
the use of Feynman diagrams almost inescapable in QFT.

An easy consequence of the preceding analysis of our Feynman
diagram structures, obtained by counting the half-lines and using
the fact that the $H$-vertices have valence at least $3$, is

\begin{lemma}
A tree-like connected Feynman diagram, which is then
necessarily of type $(I,J)$ with $\#(I)=1$,
satisfies
\be
\#({\til E})\le 2 l
\ee
where $l$ is the total number of leaves
$l\eqdef \#({\til E}_Y)+\#({\til E}_{\Br u})$.
\end{lemma}
From which one deduces by adding the above inequalities obtained
for each tree growing off the central circuit, that

\begin{lemma}
A circuit-like connected Feynman diagram, which is then
necessarily of type $(I,J)$ with $I=\emptyset$,
satisfies also
\be
\#({\til E})\le 2 \lp\#({\til E}_Y)+\#({\til E}_{\Br u})\rp
\ee
\end{lemma}
Finally by adding the inequalities for each connected component

\begin{lemma}
Any Feynman diagram, of arbitrary type $(I,J)$,
also satisfies
\be
\#(E)\le
\#({\til E})\le 2 \lp\#({\til E}_Y)+\#({\til E}_{\Br u})\rp
=2\lp\#(\pi_Y)+\#(J)\rp
\ee
\end{lemma}

Although quite trivial the above lemmas
are crucial in order to
ensure that the grading,
with respect to which the topology of the ring
$R[[Y]]$ is defined, and which is related to $Y$-vertices only, grows
with the complexity of the Feynman diagram.
This observation securing the summability and the same argument
as in Proposition 5 now entail the following.

\begin{prop}
Let the finite sets $I$ and $J$ and the assignment maps
$\ta_u$ and $\ta_{\Br u}$ be given.
Let $U$ be the element of
$R[[{\Br u},u,Y]]$ given by
\be
U=
\lp
\prod_{i\in I}
u_{\ta_u(i)}
\rp
\lp
\prod_{j\in J}
{\Br u}_{\ta_{\Br u}(j)}
\rp
\exp\lp
{\Br u}H(u)+{\Br u}Y
\rp
\ee
then the following identity holds in $R[[Y]]$,
both sides being summable
\be
\int d{\Br u}du\ e^{-{\Br u}A u}\ U=
(det\ A)^{-1}
\sum_{[E,\cE]}
\sum_{\cC}
\frac{\cA_{Fey}
\lp
E,(\cE,\cC),\ta_u,\ta_{\Br u}
\rp}
{\#Aut(E,\cE)}
\ee
where the sum is over equivalence classes of pre-Feynman diagram structures
of type $(I,J)$,
$(E,\cE)$ being an arbitrary representative of such a class.
$\cC$ is summed over contraction schemes i.e. bijective
maps $\cC:E_{\Br u}\rightarrow E_u$.
\end{prop}

\noindent
By the same proof as that of Theorem 3, one now arrives at

\begin{theor}
With the same hypothesis as in the previous proposition one has,
both sides being summable in
$R[[Y]]$,
\be
\int d{\Br u}du\ e^{-{\Br u}A u}\ U=
(det\ A)^{-1}
\sum_{[E,\cF]}
\frac{\cA_{Fey}
\lp
E,\cF,\ta_u,\ta_{\Br u}
\rp}
{\#Aut(E,\cF)}
\ee
where the sum is over equivalence classes of Feynman diagram structures
of type $(I,J)$, and $(E,\cF)$ denotes an arbitrary
class representative.
\end{theor}

We have now completely defined, in a mathematically precise fashion,
the numerator and the denominator that appear in Claim 3.
They correspond with the situation where
$(I,J)=([1],\emptyset)$ with $\ta_u(1)=i$,
and the situation where $(I,J)=(\emptyset,\emptyset)$ respectively.
The Feynman diagrams in the former situation can be called, according to
physical terminology,
{\em 1-point} diagrams.
In the latter situation they would rather be called {\em vacuum} diagrams.
Before we end this section we still have to prove the following
precise restatement of Claim 3.

\begin{theor}
The compositional inverse of $F=(F_i)_{1\le i\le n}$
satisfies in the ring $R[[Y]]$
the equation
\be
(F^{-1})_i(Y)=
\frac{\int d{\Br u}du\ e^{-{\Br u}A u}
\ u_i e^{{\Br u}H(u)+{\Br u}Y}}
{\int d{\Br u}du\ e^{-{\Br u}A u}
\ e^{{\Br u}H(u)+{\Br u}Y}}
\ee
the denominator being invertible in $R[[Y]]$.
\end{theor}

We will use the standard statistical mechanics notation $<.>$ for averages
and introduce,
given the finite sets $I$ and $J$ and their associated assignment maps
$\ta_u$ and $\ta_{\Br u}$,
the {\em unnormalized correlation function}
\bea
\lefteqn{
<
\lp
\prod_{i\in I}
u_{\ta_u(i)}
\rp
\lp
\prod_{j\in J}
{\Br u}_{\ta_{\Br u}(j)}
\rp
>_U\eqdef
} & & \nonumber \\
 & & (det\ A)
\int d{\Br u}du\ e^{-{\Br u}A u}
\lp
\prod_{i\in I}
u_{\ta_u(i)}
\rp
\lp
\prod_{j\in J}
{\Br u}_{\ta_{\Br u}(j)}
\rp
e^{{\Br u}H(u)+{\Br u}Y}
\eea
Note that
\be
det\ A=\frac{1}{\int d{\Br u}du\ e^{-{\Br u}A u} \ 1}
\ee
represents the normalization by its total weight
(in order to have a probability measure)
of the ``Gaussian measure''
$d{\Br u}du\ e^{-{\Br u}A u}$.
It is not the full ``interacting measure''
$d{\Br u}du\ e^{-{\Br u}A u +{\Br u}H(u)+{\Br u}Y}$, hence
the word ``unnormalized''.
The corresponding {\em normalized correlation function} is rather
\be
<
\lp
\prod_{i\in I}
u_{\ta_u(i)}
\rp
\lp
\prod_{j\in J}
{\Br u}_{\ta_{\Br u}(j)}
\rp
>_N\eqdef
\frac{1}{Z}
<
\lp
\prod_{i\in I}
u_{\ta_u(i)}
\rp
\lp
\prod_{j\in J}
{\Br u}_{\ta_{\Br u}(j)}
\rp
>_U
\ee
where the $Z$ is the {\em partition function}
defined by
\be
Z\eqdef
<1>_U=
(det\ A)
\int d{\Br u}du\ e^{-{\Br u}A u}
\ e^{{\Br u}H(u)+{\Br u}Y}
\ee
It is given by Theorem 4 as a sum over classes of,
not necessarily connected,
vacuum
(i.e. of type $(\emptyset,\emptyset)$)
Feynman diagram structures
\be
Z=\sum_{{[E,\cF]} \atop {{\rm type}\ (\emptyset,\emptyset)}}
\frac{\cA_{Fey}
\lp
E,\cF,\ta_u,\ta_{\Br u}
\rp}
{\#Aut(E,\cF)}
\ee
The constant term of $Z$ is given by the contribution of
the trivial diagram corresponding to $E=\emptyset$, and is equal to
$1$ (one can check that our definitions also
hold in this degenerate case).
As a result, $Z$ i.e. the denominator in Theorem 5 is invertible in
$R[[Y]]$.
One can also define the subspecie of nontrivial
connected vacuum Feynman diagrams
$\cF$ on a set $E$ by adding to Definition 12,
in the case where $I=J=\emptyset$, the condition that
$E\neq\emptyset$ and that
the associated set ${\til E}$ and digraph $G$
are such that $G$ connects ${\til E}$.
One can then define the {\em free energy}
\be
W\eqdef
\sum_{{[E,\cF]\ {\rm type}\ (\emptyset,\emptyset)}
\atop {{\rm connected}\ E\neq\emptyset}}
\frac{\cA_{Fey}
\lp
E,\cF
\rp}
{\#Aut(E,\cF)}
\ee
which is summable in $R[[Y]]$, as a part of the sum for $Z$ which
is already known to be summable.
One can also prove this directly using Lemma 3.
Note that there is no longer a need to specify the maps
$\ta_u$ and $\ta_{\Br u}$ whose graphs are empty.
Note also that
the diagrams appearing in the last equation are each made of a single
nonempty circuit-like connected component whose leaves are all $Y$-vertices.
Now it is easy to check that
\begin{prop}
\be
Z=\exp\lp W\rp
\ee
\end{prop}
Similar statements for Hurewitz or exponential generating series
are quite familiar in combinatorial theory.
It boils down to the use of the multinomial theorem,
the invariance of amplitudes by relabelling and, most importantly here,
their {\em factorization} over connected components.

\medskip
In fact, for any fixed type $(I,J)$
one can define in an analogous way,
the subspecie of connected Feynman diagram structures of type $(I,J)$,
by requiring that the digraph $G$ connects the derived set ${\til E}$.
This allows, again given the assignment maps $\ta_u$ and $\ta_{\Br u}$,
to define the {\em connected correlation functions}
\be
<
\lp
\prod_{i\in I}
u_{\ta_u(i)}
\rp
\lp
\prod_{j\in J}
{\Br u}_{\ta_{\Br u}(j)}
\rp
>_C\eqdef
\sum_{{[E,\cF]\ {\rm type}\ (I,J)}
\atop {{\rm connected}}}
\frac{\cA_{Fey}
\lp
E,\cF,\ta_u,\ta_{\Br u}
\rp}
{\#Aut(E,\cF)}
\label{connect}
\ee
These can also be called {\em cumulants} or {\em semi-invariants}
in conformity with the terminology of mathematical
statistics and probability theory.
They are also related to the so-called {\em Ursell functions}
in statistical mechanics.
Indeed, one has
\begin{theor}
\bea
\lefteqn{
<
\lp
\prod_{i\in I}
u_{\ta_u(i)}
\rp
\lp
\prod_{j\in J}
{\Br u}_{\ta_{\Br u}(j)}
\rp>_U =
} & & \nonumber \\
 & & Z\times \sum_{\pi}
\prod_{({\til I},{\til J})\in\pi}
<
\lp
\prod_{i\in {\til I}}
u_{\ta_u|_{\til I}(i)}
\rp
\lp
\prod_{j\in {\til J}}
{\Br u}_{\ta_{\Br u}|_{\til J}(j)}
\rp>_C
\label{cluster}
\eea
where the sum is over all (unordered) sets $\pi$ of pairs
$({\til I},{\til J})$
such that ${\til I}$
and ${\til J}$
are not simultaneously empty
subsets of $I$ and $J$ respectively, and such that
\be
\pi_I\eqdef
\{
{\til I}\subset I|{\til I}\neq\emptyset\ {\rm and}\ \exists
{\til J}\subset J, ({\til I},{\til J})\in\pi
\}
\ee
and
\be
\pi_J\eqdef
\{
{\til J}\subset J|{\til J}\neq\emptyset\ {\rm and}\ \exists
{\til I}\subset I, ({\til I},{\til J})\in\pi
\}
\ee
are partitions of $I$ and $J$ respectively.
\end{theor}

\noindent{\bf Proof :}
One starts from the expression given by Theorem 4
for the unnormalized correlation function
\[
<
\lp
\prod_{i\in I}
u_{\ta_u(i)}
\rp
\lp
\prod_{j\in J}
{\Br u}_{\ta_{\Br u}(j)}
\rp>_U 
\]
as a sum over classes $[E,\cF]$
of corresponding Feynman diagrams.
Given such a diagram, one divides $E$ according to the connected
components of ${\til E}$ that are determined by the
digraph $G$.
We let $E_Z\subset E$
be the union of vacuum connected components
(i.e. those which do not intersect the images of $\rh_u$ and
$\rh_{\Br u}$).
For any set of labels $F\subset E$ corresponding to a connected
component which does intersect $\rh_u(I)$
and $\rh_{\Br u}(J)$,
we consider $I_F\eqdef\rh_u^{-1}(\rh_u(I)\cap F)$ and
$J_F\eqdef\rh_{\Br u}^{-1}(\rh_{\Br u}(J)\cap F)$
and we let $\pi$ be the set of pairs $(I_F,J_F)$ obtained in this way.
The set $\pi$ satisfies the conditions stated in the theorem.
For each $({\til I},{\til J})\in\pi$, we let $E_{({\til I},{\til J})}$
be the unique component $F$ of $E$
such that ${\til I}=I_F$
and ${\til J}=J_F$.
One then canonically deduces from the Feynman diagram structure 
$\cF$ of type $(I,J)$
on $E$ an induced connected diagram structure $\cF_{({\til I},{\til J})}$
of type $({\til I},{\til J})$ on $E_{({\til I},{\til J})}$.
One also obtains in the same obvious manner a (not necessarily connected)
Feynman diagram structure $\cF_Z$ of type
$(\emptyset,\emptyset)$ on $E_Z$.
The index assignment maps for a pair $({\til I},{\til J})\in\pi$
are defined from $\ta_u$ and $\ta_{\Br u}$ by restriction from
$I$ to ${\til I}$ and from $J$ to ${\til J}$
respectively. All one has to do  in proving the equality (\ref{cluster})
is to notice that one can replace the global sum over
$[E,\cF]$
by the sum over the set $\pi$ and independent sums on the
classes $[E_{({\til I},{\til J})},\cF_{({\til I},{\til J})}]$
for $({\til I},{\til J})\in\pi$ and the class $[E_Z,\cF_Z]$,
the amplitudes being factorized over connected components
and also the symmetry factors.
Indeed one has a canonical group isomorphism
\be
Aut(E,\cF)
\simeq
Aut(E_Z,\cF_Z)\times\prod_{({\til I},{\til J})\in\pi}
Aut\lp
E_{({\til I},{\til J})},\cF_{({\til I},{\til J})}
\rp
\ee
\endproof

\noindent
An immediate consequence is that

\noindent{\bf Corollary 1 :}
\be
\frac{\int d{\Br u}du\ e^{-{\Br u}A u}
\ u_i e^{{\Br u}H(u)+{\Br u}Y}}
{\int d{\Br u}du\ e^{-{\Br u}A u}
\ e^{{\Br u}H(u)+{\Br u}Y}}
= <u_i>_C
\ee
which is a sum over connected Feynman diagrams of type $([1],\emptyset)$,
which must be tree-like with leaves excusively made of $Y$-vertices.

\medskip
\noindent{\bf Proof of Theorem 5 :}
One starts from
\be
<u_i>_C=
\sum_{{[E,\cF]\ {\rm type}\ ([1],\emptyset)}
\atop {\rm connected}}
\frac{\cA_{Fey}
\lp
E,\cF,\ta_u,\ta_{\Br u}
\rp}
{\#Aut(E,\cF)}
\ee
with $\ta_u(1)=i$ and
$\ta_{\Br u}$ empty.
In the previous sum one distinguishes the simplest
term corresponding to the diagram class
\be
\figput{dessin10}
\ee
for which $\#(E)=2$, $\#Aut(E,\cF)=1$
and the amplitude is given by
\be
\cA_{Fey}
\lp
E,\cF,\ta_u,\ta_{\Br u}
\rp
=\sum_{j=1}^n
(A^{-1})_{ij} Y_j
\ee
It corresponds to the linear term of the formal inverse
$(F^{-1})_i(Y)$.
Let $\Ga$ denote the sum of the remaining terms
for which $\pi_H\neq\emptyset$. For such a term,
there is a distinguished $H$-vertex $B_0\in\pi_H$, which is closest
to the root in $E_u\cap E_{ext}$, and with $p\ge 2$  attached
tree-like connected
Feynman diagram structures of type $([1],\emptyset)$, we denote
by $(E_1,\cF_1),\ldots,(E_p,\cF_p)$.
Let $\cC$ be the set of isomorphism classes
$[E,\cF]$ of  connected
Feynman diagram structures of type $([1],\emptyset)$.
There is a bijective correspondance between classes $[E,\cF]$
appearing in $\Ga$
and finitely supported families $(m_c)_{c\in\cC}$ of integers
$m_c\in\NN$ such that
$\sum_{c\in\cC} m_c\ge 2$, defined by letting
$m_c$ count the number of indices $q$ such that $(E_q,\cF_q)$
belongs to the class $c$.
Besides the cardinal of $Aut(E,\cF)$
is completely determined by
$(m_c)_{c\in\cC}$.
So is the  amplitude of $(E,\cF)$
which we denote then by $\cA((m_c)_{c\in\cC})$.
One has trivially
\be
\#Aut(E,\cF)=
\lp
\prod_{c\in\cC}
\#Aut(c)^{m_c}
\rp
\lp
\prod_{c\in\cC} m_c!
\rp
\ee
since an isomorphism of the big tree $(E,\cF)$
can
operate inside each of the branches $(E_1,\cF_1),\ldots,(E_p,\cF_p)$
and can also exchange isomorphic branches.
One can therefore write
\be
\Ga=\sum_{p\ge 2}\sum_{{(m_c)_{c\in\cC}}\atop{\sum_{c\in\cC} m_c=p}}
\frac{\cA\lp(m_c)_{c\in\cC}\rp}
{\lp
\prod_{c\in\cC}
\#Aut(c)^{m_c}
\rp
\lp
\prod_{c\in\cC} m_c!
\rp}
\ee
which by the multinomial theorem amounts to the same thing as summing over
sequences $([E_1,\cF_1],\ldots,[E_p,\cF_p])$ of elements of $\cC$.
Thus
\be
\Ga=\sum_{p\ge 2}\frac{1}{p!}
\sum_{([E_1,\cF_1],\ldots,[E_p,\cF_p])}
\frac{\cA\lp(m_c)_{c\in\cC}\rp}
{\lp
\prod_{c\in\cC}
\#Aut(c)^{m_c}
\rp}
\ee
where $(m_c)_{c\in\cC}$ is the family of multiplicities defined
$([E_1,\cF_1],\ldots,[E_p,\cF_p])$.
This becomes
\be
\Ga=\sum_{p\ge 2}\frac{1}{p!}
\sum_{([E_1,\cF_1],\ldots,[E_p,\cF_p])}
\sum_{j,j_1,\ldots,j_p=1}^n
(A^{-1})_{ij}
H^{[p]}_{j,j_1\ldots j_p}
\prod_{q=1}^p
\lp
\frac{\cA(E_q,\cF_q,j_q)}{\#Aut(E_q,\cF_q)}
\rp
\ee
where the index $j_q$ defines the $\ta_u$ assignment map
for the subdiagram $(E_q,\cF_q)$ also of type $([1],\emptyset)$.
Noting that by definition
$H^{[p]}_{j,j_1\ldots j_p}=-F^{[p]}_{j,j_1\ldots j_p}$
and using (\ref{connect}), the previous
expression recombines into
\be
\Ga=-\sum_{p\ge 2}\frac{1}{p!}
\sum_{j,j_1,\ldots,j_p=1}^n
(A^{-1})_{ij}
F^{[p]}_{j,j_1\ldots j_p}
<u_{j_1}>_C\ldots<u_{j_p}>_C
\ee
therefore
\be
<u_i>_C
=
\sum_{j=1}^n
(A^{-1})_{ij} Y_j
-\sum_{p\ge 2}\frac{1}{p!}
\sum_{j,j_1,\ldots,j_p=1}^n
(A^{-1})_{ij}
F^{[p]}_{j,j_1\ldots j_p}
<u_{j_1}>_C\ldots<u_{j_p}>_C
\ee
Multiplying on the left by the matrix $A=(F^{[1]}_{i,j})_{1\le i,j\le n}$
and transposing the sum over $p$ gives
\be
F_i(<u>_C)=Y_i
\ee
which shows that $<u_i>_C\in R[[Y]]$
is the $i$-th component of the right compositional inverse,
that is simply the inverse, of $F=(F_i)_{1\le i\le n}$,
which concludes our proof.
\endproof

\subsection{Lagrange-Good inversion}

In order to avoid lengthy
repetitions of the previous arguments,
we will be rather brief, in this section,
and only detail the new
ingredients needed.
We work in the ring $R[[X_1,\ldots,X_n]]$.
We suppose that we have $n$ power series $(G_i)_{1\le i\le n}$
in $n$ variables defined by their tensor elements
$G^{[p]}_{i,j_1\ldots j_p}$
with $p\ge 0$.
We define as before the unnormalized correlation
functions
\be
<
\lp
\prod_{i\in I}
u_{\ta_u(i)}
\rp
\lp
\prod_{j\in J}
{\Br u}_{\ta_{\Br u}(j)}
\rp>_U =
\int d{\Br u}du\ e^{-{\Br u}u}
\lp
\prod_{i\in I}
u_{\ta_u(i)}
\rp
\lp
\prod_{j\in J}
{\Br u}_{\ta_{\Br u}(j)}
\rp
e^{{\Br u}XG(u)}
\ee
in the ring $R[[X_1,\ldots,X_n]]$,
by extending formal Gaussian integration with the identity matrix
as a covariance,
from monomials
in the $u$'s and ${\Br u}$'s to elements of $R[[{\Br u},u,X]]$,
whenever the summation
(over the multiindices defining the monomials)
converges in $R[[X]]$.
One needs almost the same definitions of pre-Feynman and Feynman
diagram structures as in section III.2 except that one
has only one type of vertices we call $XG$-vertices,
corresponding to a partition
$\pi_{XG}$ of $E_{int}$.
A block $B\in\pi_{XG}$
must have exactly one element in $E_{\Br u}$ but any number of elements
of $E_u$ is allowed this time, even zero (corresponding to the
tree leaves).
The contribution of such a $XG$-vertex in the amplitude of a
Feynman diagram is
\be
\figplace{dessin11}{0 in}{-0.45 in}
=X_i G^{[p]}_{i,j_1\ldots j_p}
\label{FruleXG}
\ee
A contraction line corresponds to a factor
\be
\figplace{dessin12}{0 in}{-0.065 in}=\de_{ij}
\label{Frulede}
\ee
in the amplitude of a Feynman
graph.
Apart from this small difference, the treatment is exactly the same
as in section III.2.
One has an analog of Corollary 1
\begin{prop}
\be
\frac{\int d{\Br u}du\ e^{-{\Br u} u}
\ u_i e^{{\Br u}XG(u)}}
{\int d{\Br u}du\ e^{-{\Br u} u}
\ e^{{\Br u}XG(u)}}
= <u_i>_C
\eqdef
\sum_{{[E,\cF]\ {\rm type}\ ([1],\emptyset)}
\atop {\rm connected}}
\frac{\cA_{Fey}
\lp
E,\cF,\ta_u,\ta_{\Br u}
\rp}
{\#Aut(E,\cF)}
\label{lag1pt}
\ee
with $\ta_u(1)=i$ and $\ta_{\Br u}=\emptyset$,
and the amplitude $\cA_{Fey}(E,\cF,\ta_u,\ta_{\Br u})$ is defined
using the {\em Feynman rules} (\ref{FruleXG}) and (\ref{Frulede}).
\end{prop}

For example the amplitude 
\bea
\lefteqn{
\cA_{Fey}
\lp
E,\cF,\ta_u,\ta_{\Br u}
\rp
=\sum_{\al_1,\ldots,\al_{10}}^n
X_i G^{[3]}_{i,\al_1\al_2\al_3}
} & & \nonumber \\
 & & \times
\lp
X_{\al_1} G^{[1]}_{\al_1,\al_4}
X_{\al_4} G^{[3]}_{\al_4,\al_5\al_6\al_7}
X_{\al_5} G^{[0]}_{\al_5}
X_{\al_6} G^{[0]}_{\al_6}
X_{\al_7} G^{[0]}_{\al_7}
\rp
\nonumber \\
 & & \times
\lp
X_{\al_2} G^{[0]}_{\al_2}
\rp
\lp
X_{\al_3} G^{[3]}_{\al_3,\al_8\al_9\al_{10}}
X_{\al_8} G^{[0]}_{\al_8}
X_{\al_9} G^{[0]}_{\al_9}
X_{\al_{10}} G^{[0]}_{\al_{10}}
\rp
\eea
is assigned to the Feynman diagram

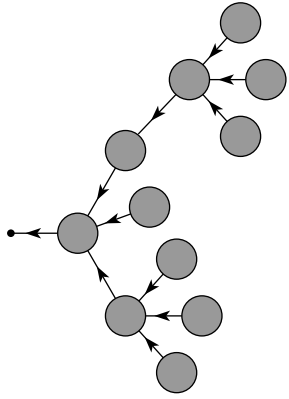
\be
\figput{dessin13}
\ee
whose automorphism group
has cardinality
$\#Aut(E,\cF)=3!\times 3!$.

The convergence of the Feynman diagram expansions in the ring $R[[X]]$
is ensured by the fact
that each vertex (and not only the leaves
like in the previous section)
increases the grading by one unit.
By repeating the same arguments as in the proof of Theorem 5,
consisting in identifying the nearest $XG$-vertex
to the root and summing over the sub-trees
that are attached to it, it is immediate that the series
\be
F_i(X)\eqdef <u_i>_C
\label{verifeq}
\ee
is a solution of the implicit equations
\be
F_i(X)=X_i G_i(F(X))\ \ {\rm for}\ \ 1\le i\le n
\ee
This gives a rigorous restatement of Claim 4 in section II.3.
One also has an analog of Proposition 9
saying that
\be
Z\eqdef
\int d{\Br u}du\ e^{-{\Br u} u}
e^{{\Br u}XG(u)}
=\exp\lp W\rp
\label{linked}
\ee
with
\be
W\eqdef
\sum_{{[E,\cF]\ {\rm type}\ (\emptyset,\emptyset)}
\atop {{\rm connected}\ E\neq\emptyset}}
\frac{\cA_{Fey}
\lp
E,\cF
\rp}
{\#Aut(E,\cF)}
\label{doublev}
\ee
where the sum is over equivalence classes of
nonempty connected vacuum Feynman diagrams.
A closer look at these diagrams will allow us to prove the
following.
\begin{theor}
Using the notations of section II.3
\be
Z=\frac{1}{det\lp·I-X\partial G(F)\rp}
\ee
\end{theor}
The proof of this theorem depends on the following result
which deserves to be stated as an independent theorem.
The argument
must be familiar to the practitioner
of combinatorial species but we could not find it stated explicitely
as we need it,
in the literature.
We cannot resist calling it ``the principle of variation of ambiguity''
and it states a kind  of {\em functoriality}
of Feynman diagrammatic perturbation series.
Ambiguity refers to the ``degree of resolution'' of
the combinatorial description that we (``the observer'') chose
and which is like a combinatorialist's ``choice of coordinates''.

\begin{theor}
Let $\cM$ and $\cN$ be two
combinatorial species in the sense of Joyal~\cite{Joyal}, related
by a natural transformation (or a morphism of functors)
$\rh$.
That is for every
finite set
$E$ we have a (not necessarily bijective) map $\rh_E:\cM(E)\rightarrow
\cN(E)$, such that for any bijection $\si:E\rightarrow F$
between finite sets $E$ and $F$, one has
$\rh_F\circ\cM(\si)=\cN(\si)\circ\rh_E$.
Suppose we have defined for every pair $(E,M)$, consisting of a finite set
$E$ and a structure $M\in\cM(E)$ of type $\cM$ on $E$,
an amplitude $\cA(E,M)$ taking values in a formal power series
ring $R[[V]]$, where $V$ denotes any collection of indeterminates
and the ground ring $R$ contains $\QQ$.
Assume that $\cA(E,M)$
is constant over equivalence classes, denoted by $[E,M]$,
of pairs $(E,M)$ for the relation
$(E,M)\sim(E',M')$ if and only if there exists a bijection
$\si:E\rightarrow E'$
with $\cM(\si)(M)=M'$.

The conclusion of the theorem is that if the left-hand side of
\be
\sum_{[E,M]}
\frac{\cA(E,M)}
{\#Aut(E,M)}
=
\sum_{[E,N]}
\frac{1}
{\#Aut(E,N)}
\sum_{{M\in\cM(E)}\atop{\rh_E(M)=N}}
\cA(E,M)
\label{variation}
\ee
converges in $R[[V]]$, then so does the left-hand
side and the equality holds.
Note that in the left-hand side one sums over classes for
the specie $\cM$, while in the right-hand side one sums over classes
for the specie $\cN$.
\end{theor}

\noindent{\bf Proof :}
Note that by the equivariance of the transformation $\rh$, the
expression
\[
\frac{1}
{\#Aut(E,N)}
\sum_{{M\in\cM(E)}\atop{\rh_E(M)=N}}
\cA(E,M)
\]
is independent of the pair $(E,N)$ in a given class $[E,N]$
for the specie $\cN$.
Note also that there is no set-theoretic difficulty in speaking of
``the set of all equivalence classes'' for a specie $\cM$.
Indeed such a set can be easily constructed as a quotient of the
disjoint union of the denumerable family of finite sets
$(\cM([k]))_{k\in\NN}$.
Therefore the families of elements of $R[[V]]$
to be summed in both sides of (\ref{variation})
are well-defined.
Let us first show that the summability of the left-hand side
implies that of the right-hand side.
One can define a map ${\Br \rh}$ from the set of equivalence classes
of $\cM$ to that
of $\cN$ by
\be
{\Br \rh}([E,M])\eqdef
[E,\rh_E(M)]
\ee
Let $V^\al$ be a monomial in $R[[V]]$ and denote by
$[V^\al]\Om$ the coefficient of $V^\al$ in a power series
$\Om\in R[[V]]$.
If
\be
[V^\al]
\lp
\frac{1}
{\#Aut(E,N)}
\sum_{{M\in\cM(E)}\atop{\rh_E(M)=N}}
\cA(E,M)
\rp
\neq 0
\label{nonnul}
\ee
then $[E,N]$
is the image by ${\Br \rh}$
of a $[E,M]$
such that $[V^\al]\cA(E,M)\neq 0$.
If the left-hand side converges, there
are finitely many
such $[E,M]$'s
and since ${\Br \rh}$ is a finite-to-one map, there
are finitely $[E,N]$'s
such that (\ref{nonnul}) is true.

To prove the equality in (\ref{variation}),
one simply needs to check that for any class
$[E,N]$ for the specie $\cN$
the following equality, involving only finite sums, holds:
\be
\sum_{{[E,M]}\atop{{\Br \rh}([E,M])=
[E,N]}}
\frac{\cA(E,M)}
{\#Aut(E,M)}
=
\frac{1}
{\#Aut(E,N)}
\sum_{{M\in\cM(E)}\atop{\rh_E(M)=N}}
\cA(E,M)
\label{varfin}
\ee
First fix a representative $(E,N)$ of the concerned $\cN$-class.
Let $\pi_\cN$ be the partition of $\cN(E)$
into equivalence classes for the relation
$N_1\sim N_2$
defined by the existence of
a bijection $\si:E\rightarrow E$
such that $\cN(\si)(N_1)=N_2$.
Let $\pi_\cM$ be the analogous partition of $\cM(E)$.
Let $\pi_\rh$ be the partition of $\cM(E)$ defined by the nonempty
inverse images by $\rh_E$ of elements of $\cN(E)$.
It is clear by functoriality of $\rh$ that $\pi_\cM$ is finer
than $\pi_\rh$.
Let ${\Br N}$
be the block of $\pi_\cN$ containing $N$.
One easily check
\be
\sum_{{[E,M]}\atop{{\Br \rh}([E,M])=
[E,N]}}
\frac{\cA(E,M)}
{\#Aut(E,M)}
=
\sum_{{B\in\pi_\cM}\atop{\rh_E(B)\subset {\Br N}}}
\frac{1}
{\#(E)!}
\frac{\#(E)!}{\#Aut(E,M)}
\cA(E,M)
\label{interm}
\ee
where $M$ designates any element of $B$.
Indeed every class $[E,M]$
with ${\Br \rh}([E,M])=[E,N]$
corresponds to a $B\in\pi_\cM$ sent
by $\rh_E$ into ${\Br N}$.
One also has
\be
\frac{\#(E)}{\#Aut(E,M)}=\#(B)
\ee
since $B$ is the orbit of any
$M\in B$ for the action of $\GS(E)$, the group of
permutations of $E$, on the set $\cM(E)$.
Therefore the right hand side of (\ref{interm})
becomes
\be
\frac{1}
{\#(E)!}
\sum_{{M\in\cM(E)}\atop{\rh_E(M)\in{\Br N}}}
\cA(E,M)
=
\frac{1}
{\#(E)!}
\sum_{N'\in {\Br N}}
\sum_{{M\in\cM(E)}\atop{\rh_E(M)=N'}}
\cA(E,M)
\ee
Now again by functoriality of $\rh$
\[
\sum_{{M\in\cM(E)}\atop{\rh_E(M)=N'}}
\cA(E,M)
\]
does not depend on
$N'$ in ${\Br N}$ the latter of which is
the orbit of $N$ under the action of 
$\GS(E)$ on $\cN(E)$, and therefore has
cardinality
\[
\frac{\#(E)!}{\#Aut(E,N)}
\]
Thus
\be
\frac{1}
{\#(E)!}
\sum_{{M\in\cM(E)}\atop{\rh_E(M)\in{\Br N}}}
\cA(E,M)
=
\frac{1}{\#Aut(E,N)}
\sum_{{M\in\cM(E)}\atop{\rh_E(M)=N}}
\cA(E,M)
\ee
from which (\ref{varfin}) and the proof of the theorem follow.
\endproof

\noindent{\bf Remark :}
We have already used this principle in two particular instances:
\begin{itemize}
\item
In Theorem 7, where $\cM$ was the species of Feynman diagrams,
and $\cN$ that of pre-Feynman diagrams. The transformation $\rh$
amounted to forgetting
the contraction scheme.
\item
In (\ref{Hurewitz}), where $\cM$ was the specie of Feynman diagrams,
and $\cN$ was the vacuous specie
($\#(\cN(E))=1$ for any finite $E$).
Applying $\rh$ meant to forget everything
except the cardinality
of $E$.
\end{itemize}

\medskip
\noindent{\bf Proof of Theorem 7 :}
We start from the expression (\ref{doublev}) for $W$ that we rewrite,
following
the notation of Theorem 8, as
\be
W=\sum_{[E,N]}
\frac{\cA(E,N)}
{\#Aut(E,N)}
\ee
Here the species $\cN$ is that of nontrivial connected Feynman diagram
structures of type $(\emptyset,\emptyset)$.
The amplitude is the one defined by the Feynman rules
(\ref{FruleXG}) and (\ref{Frulede}).
We now introduce a new specie $\cM$
as follows.
For any finite set $E$, we call an $\cM$-structure on $E$,
any couple $(N,\cO)$ where $N\in\cN(E)$
and $\cO$ consists
of a total ordering
$B_1<\ldots< B_p$ of the $XG$-vertices appearing in the central
circuit of $N$,
and of a total ordering $x_1^q<\ldots<x_{k_q}^q$
of the elements of $B_q\cap E_u$ for each
$q$, $1\le q\le p$.
We require that the order
$B_1<\ldots< B_p$ be compatible with the orientation of the circuit, i.e.
$B_1,\ldots,B_p$ is
the sequence of vertices obtained by following the orientation
of the contraction lines, along the circuit, starting from $B_1$.
Note that it is possible that
some $k_q\eqdef \#(B_q\cap E_u)$
be zero.
However one allways has $p\ge 1$.
Transport of structure for $\cM$ is defined in the obvious
covariant way.
The morphism of functors $\rh$ is defined by
$\rh_E(N,\cO)=N$
for any $(N,\cO)\in\cM(E)$.
We also define the amplitude for an $\cM$-structure
$M=(N,\cO)$, keeping the previous
notations, by
\be
\cA(E,M)
\eqdef
\frac{\cA(E,N)}{p. k_1!\ldots k_p!}
\ee
Now, after the trivial check that the left hand side
of the following equality converges in $R[[X]]$,
Theorem 8
implies that
\be
\sum_{[E,M]}
\frac{\cA(E,M)}
{\#Aut(E,M)}
=
\sum_{[E,N]}
\frac{1}
{\#Aut(E,N)}
\sum_{{M\in\cM(E)}\atop{\rh_E(M)=N}}
\cA(E,M)
\ee
But
\be
\sum_{{M\in\cM(E)}\atop{\rh_E(M)=N}}
\cA(E,M)
=
\cA(E,N)
\ee
Indeed, the {\em product} $p.k_1!\ldots k_p!$ does not depend on
$\cO$,
besides
it is equal to the number of these
possible orderings $\cO$.
One can write as a result
\be
W=\sum_{[E,M]}
\frac{\cA(E,\rh_E(M))}{p.k_1!\ldots k_p!}
\frac{1}{\#Aut(E,M)}
\ee
The point is that the automorphism group of a pair
$(E,M)$ is much more manageble since
the central circuit
has been completely
{\em rigidified}, i.e.
all the elements of $E$ that belong to an $XG$-vertex
along the circuit are fixed by automorphisms of $(E,M)$.
Indeed, for any $q$, $1\le q\le p$ and any $\nu$,
$1\le \nu\le k_q$,
$x_\nu^q\in B_q\cap E_u$ is the
new root of a tree-like connected Feynman diagram
structure $\cF_\nu^q$ of type $([1],\emptyset)$ on a subset
$E_\nu^q$
of $E$.
The corresponding $\rh_u$ map has $\{x_\nu^q\}$
as an image.
An automorphism of $(E,M)$ has to restrict
inside
$E_\nu^q$
to an automorphism of $\cF_\nu^q$.
Therefore
\be
\#Aut(E,M)
=
\prod_{q=1}^p
\lp
\prod_{\nu=1}^{k_q}
\#Aut(E_\nu^q,\cF_\nu^q)
\rp
\ee
Besides the amplitude $\cA(E,N)$,
with $N=\rh_E(M)$
is given by
\be
\cA(E,N)
=
\sum_I
\cL_I
\prod_{q=1}^p
\lp
\prod_{\nu=1}^{k_q}
\cA(E_\nu^q,\cF_\nu^q,i_\nu^q)
\rp
\ee
where the sum is over families $I=(i_\nu^q)_{1\le q\le p, 1\le \nu\le k_q}$
of indices
in $[n]$.
$\cA(E_\nu^q,\cF_\nu^q,i_\nu^q)$ is the amplitude
of the Feynman diagram structure $\cF_\nu^q$ of
type $([1],\emptyset)$
on $E_\nu^q$
with respect to the index assignment
map with value $i_\nu^q$.
Finally $\cL_I$ is the contribution of the {\em amputated} circuit
\be
\cL_I\eqdef
\sum_{j_1,\ldots,j_p=1}^n
\prod_{q=1}^p
\lp
X_{j_q}
G^{[k_q+1]}_{j_q,j_{q+1}i_1^q\ldots i_{k_q}^q}
\rp
\ee
with the convention that $j_{p+1}\eqdef j_1$.
Note also that classes $[E,M]$ are in bijective correspondance with
families $([E_\nu^q,\cF_\nu^q])_{1\le q\le p, 1\le \nu\le k_q}$
of classes of connected Feynman diagram structures
of type $([1],\emptyset)$
where all values of $p\ge 1$ and $k_q\ge 0$, for
$1\le q\le p$, are allowed.
The previous observation,
equations (\ref{lag1pt}) and (\ref{verifeq}),
and the
expression of a derivative $\partial_j G_i$
in tensorial notation is enough
to show that
\be
W=
\sum_{p\ge 1}
\frac{1}{p}
tr
{\left[
X\partial G\lp F(X)\rp
\right]}^p
\ee
and, as a result of Jacobi's identity and equation (\ref{linked})
\be
Z=\exp\lp W\rp=
\frac{1}{det\lp
I-X\partial G(F)
\rp}
\ee 
\endproof

Note that we have an analog of Theorem 6 whose statement
and proof are the same in the present
context.
As a consequence one has
\begin{theor}
For any monomial
$\Om(F)=F_1^{\al_1}\ldots F_n^{\al_n}$,
the following identity in
$R[[X]]$,
both sides being well-defined, holds
\be
\Om(F)
\times
\frac{1}{det\lp
I-X\partial G(F)
\rp}
=
\int d{\Br u}du\ e^{-{\Br u} u}
\ \Om(u) e^{{\Br u}XG(u)}
\ee
\end{theor}

\noindent
It easily checked that one can expand the $e^{{\Br u}XG(u)}$
and take out
the sum to get, in the ring $R[[X]]$:
\be
\int d{\Br u}du\ e^{-{\Br u} u}
\ \Om(u) e^{{\Br u}XG(u)}
=
\sum_{\al\in\NN^n}
\frac{X^\al}{\al!}
\int d{\Br u}du\ e^{-{\Br u} u}
\ {\Br u}^\al
\Om(u)
G(u)^\al
\ee
Note that
\be
\int d{\Br u}du\ e^{-{\Br u} u}
\ {\Br u}^\al
\Om(u)
G(u)^\al
\in
R
\ee
and can be computed, by going back to the definition
of formal Gaussian integration with covariance matrix given by
the identity matrix, as
\[
{\left.
{\lp
\frac{\partial}
{\partial u}
\rp}^\al
\right|}_{u=0}
\left[
\Om(u)
G(u)^\al
\right]
\]
This concludes our derivation of the implicit form of the
multivariable Lagrange-Good inversion.

\section{Comments}

\noindent{\bf 1)}
By now, it should be clear to the reader that we have only
scratched the tip
of the iceberg.
In QFT, there are basically four categories of fields (i.e. types
of variables on which one can define a formal Gaussian integration scheme).
This division is strangely reminiscent of the distinction
between the main families of classical groups. We indeed have :

\begin{itemize}
\item
Complex Bosonic fields :
The variables commute and come with an involution exchanging them in pairs.
This is the situation we covered here.
Expansions involve digraphs, and Wick's theorem is in terms
of {\em permanents}.
\item
Complex Fermionic fields :
The variables anti-commute and also come with an involution.
Graphs are directed but usually involve
an extra $-1$ factor per circuit.
Wick's theorem uses {\em determinants}.
In many respects, Fermionic integration intuitively behaves like
Bosonic integration in a
``negative dimensional space'', whatever that means.
\item
Real Bosonic fields :
The variables commute and no involution on them
is given.
The graphs are undirected.
Wick's theorem is in term of {\em hafnians},
i.e. sums are over perfect matchings instead of permutations.
The covariance matrices must be symmetric.
\item
Real Fermionic fields :
The variables anti-commute.
No involution is, at least beforehand, given.
Covariance matrices must be skew-symmetric, therefore
graphs have to be, somewhat
artificially, oriented to avoid sign
ambiguities in their amplitudes.
Wick's theorem involves {\em Pfaffians}.
\end{itemize}
Clearly, a similar approach to ours, using
combinatorial species, can be developped for all four
types of fields; although
one has to be careful with Fermions.
For instance, we do not know if one can make sense
of situations where vertices have an odd number of half-lines or,
in the complex case, unequal numbers of incoming and outgoing half-lines.
The case of ribbon graphs (see~\cite{Fiorenza} for instance)
is covered by the above tentative classification.
The GUE random matrix ensemble, for example, belongs
to the complex Bosonic case, while the GOE falls in the
real Bosonic case.

\noindent{\bf 2)}
Rules 1 and 2 of our symbolic calculus are rather tautological
on the diagrammatic side; but
Rule 3 can be
understood as a set of {\em combinatorial conjectures}.
Indeed we only proved the correctness of the change of variable
formula in a few special cases.
It would be a valuable task to explore the extent of its validity.
Since determinants are involved in the Jacobian factor, and thus
possibly Fermions anyway,
it might be a good idea
to directly attempt a Feynman diagrammatic
statement and proof of its
supersymmetric
generalization:
the Berezin change of variable formula.
For someone unfamiliar with this beautiful identity, we recommend
consulting :
the appendix A of~\cite{SjostrandW} which is a very clear and concise
``formulaire raisonn\'e'' of supersymmetry; then the second chapter
of~\cite{Efetov}
to see some examples of calculations and get some
exposure to
the difficulties
due to boundary terms (which however should not
intervene for what we have in mind
since, to have a Feynman diagram expansion, one needs to integrate
over the whole Bosonic space in the presence of a Gaussian weight);
and finally~\cite{Berezin} for a thorough exposition.

\noindent{\bf 3)}
Another oddity of the complex Bosonic situation we treated
here is that fields or variables come in pairs
${\Br u},u$. As our starting point was Theorem 1,
we have thought of ${\Br u}$ and $u$
as {\em complex conjugate} of one another, and have designed our notation
accordingly.
However it seems, with respect to the
change of variable formula, that
${\Br u}$ and $u$ can be manipulated as {\em independent} variables.
In fact, Rule 2 which is
a kind of Fourier representation of the Dirac delta function,
rather suggests one think of ${\Br u}$ and $u$
as {\em Fourier-dual} 
variables.
Indeed, one can derive the Gurjar-Abhyankar
formula for the formal inverse
of a system of power series in the latter spirit
by a moderate use of the theory of
pseudodifferential and Fourier integral operators
(see Exercise 3.2 in~\cite{GrigisS}).

\noindent{\bf 4)}
There is a definite and quite strange mixture of mathematics with
{\em metamathematics}
in Feynman diagrammatic sums.
As we mentioned earlier,
to describe these expansions in a mathematically precise way,
one has to define a ``programming language''
with its syntactic rules.
The sum over diagrams is in fact a sum over ``programs''
of the ``number'' (i.e. the amplitude)
such a program is meant to compute.
Some might think that this is too far-fetched an analogy,
and that basic graph theory is enough to accomodate QFT.
This is not quite correct.
In constructive field theory, the most powerful tools
are the so-called {\em phase-cell} or {\em multiscale cluster
expansions} (see~\cite{Balaban,MagnenRS,Rivasseau}
for the current state-of-the-art).
These are a kind of {\em smart} perturbation theory
designed to avoid all divergences that
appear in the naive perturbative QFT.
They make critical use of two extra ingredients:
the Heisenberg uncertainty principle (``cluster expansion'' refers
to the implementation of this idea), and
the Wilsonian renormalization group (to which ``multiscale'' refers).
We can assure the reader that
the, quite formidable,
combinatorial structures that appear in the explicit form of these
expansions~\cite{Abdesselam3},
look much more like ``programs'' than graphs.
Had we known of the theory of species at the time, we
would have written
what we called ``Mayer configurations''
(in chapter 4 of~\cite{Abdesselam3})
in this most convenient language.
Let us finish, by saying that this
intrusion of metamathematics in a problem of mathematical analysis
and also its somewhat reckless
treatment in the physical 
literature, rather than the lack of concepts (of which
the genius of K. Wilson has provided an ample supply)
is the main reason delaying
the entry of what we called the ``grammar'' of QFT
into mainstream mathematics.
Because of this, 
we venture to say that, maybe, it is time for professionals to step in:
combinatorialists, computer scientists and,
why not, mathematical logicians!

\end{document}